\documentclass[12pt]{article}
\usepackage{theorem}
\usepackage{amssymb}
\usepackage{graphicx}
\usepackage[mathscr]{eucal}
\usepackage{amsbsy}
\usepackage{amsmath}
\newtheorem{prop}{}[section]

{\theorembodyfont{\upshape} \newtheorem{rema}[prop]{}}
\newcommand{\boma}[1]{{\mbox{\boldmath $#1$} }}

\begin{document}
\newcommand{\uper}[1]{\stackrel{\barray{c} {~} \\ \mbox{\footnotesize{#1}}\farray}{\longrightarrow} }
\newcommand{\nop}[1]{ \|#1\|_{\piu} }
\newcommand{\no}[1]{ \|#1\| }
\newcommand{\nom}[1]{ \|#1\|_{\meno} }
\newcommand{\UU}[1]{e^{#1 \AA}}
\newcommand{\UD}[1]{e^{#1 \Delta}}
\newcommand{\bb}[1]{\mathbb{{#1}}}
\newcommand{\HO}[1]{\bb{H}^{{#1}}}
\newcommand{\Hz}[1]{\bb{H}^{{#1}}_{\zz}}
\newcommand{\Hs}[1]{\bb{H}^{{#1}}_{\ss}}
\newcommand{\Hg}[1]{\bb{H}^{{#1}}_{\gg}}
\newcommand{\HM}[1]{\bb{H}^{{#1}}_{\so}}
\def\ef{\psi}
\def\fun{\mathcal{F}}
\def\fun{{\tt f}}
\def\tvainf{\vspace{-0.4cm} \barray{ccc} \vspace{-0,1cm}{~}
\\ \vspace{-0.2cm} \longrightarrow \\ \vspace{-0.2cm} \scriptstyle{T \vain + \infty} \farray}
\def\De{F}
\def\er{\epsilon}
\def\erd{\er_0}
\def\Tn{T_{\star}}
\def\Tc{T_{\tt{c}}}
\def\Tb{T_{\tt{b}}}
\def\Tl{\mathscr{T}}
\def\Tm{T}
\def\Ta{T_{\tt{a}}}
\def\ua{u_{\tt{a}}}
\def\Tg{T_{G}}
\def\Tgg{T_{I}}
\def\Tw{T_{w}}
\def\Ts{T_{\Ss}}
\def\Tr{\Tl}
\def\Sp{\Ss'}
\def\Tsp{T_{\Sp}}
\def\vsm{\vspace{-0.1cm}\noindent}
\def\comple{\scriptscriptstyle{\complessi}}
\def\ug{u_G}
\def\nume{0.407}
\def\numerob{0.00724}
\def\deln{7/10}
\def\delnn{\dd{7 \over 10}}
\def\e{c}
\def\p{p}
\def\z{z}
\def\symd{{\mathfrak S}_d}
\def\del{\omega}
\def\Del{\delta}
\def\Di{\Delta}
\def\Ss{{\mathscr{S}}}
\def\Ww{{\mathscr{W}}}
\def\mmu{\hat{\mu}}
\def\rot{\mbox{rot}\,}
\def\curl{\mbox{curl}\,}
\def\Mm{\mathscr M}
\def\XS{\boma{x}}
\def\TS{\boma{t}}
\def\Lam{\boma{\eta}}
\def\DS{\boma{\rho}}
\def\KS{\boma{k}}
\def\LS{\boma{\lambda}}
\def\PR{\boma{p}}
\def\VS{\boma{v}}
\def\ski{\! \! \! \! \! \! \! \! \! \! \! \! \! \!}
\def\h{L}
\def\EM{M}
\def\EMP{M'}
\def\R{R}
\def\Rr{{\mathscr{R}}}
\def\Zz{{\mathscr{Z}}}
\def\E{E}
\def\FFf{\mathscr{F}}
\def\A{F}
\def\Xim{\Xi_{\meno}}
\def\Ximn{\Xi_{n-1}}
\def\lan{\lambda}
\def\om{\omega}
\def\Om{\Omega}
\def\Sim{\Sigm}
\def\Sip{\Delta \Sigm}
\def\Sigm{{\mathscr{S}}}
\def\Ki{{\mathscr{K}}}
\def\Hi{{\mathscr{H}}}
\def\zz{{\scriptscriptstyle{0}}}
\def\ss{{\scriptscriptstyle{\Sigma}}}
\def\gg{{\scriptscriptstyle{\Gamma}}}
\def\so{\ss \zz}
\def\Dv{\bb{\DD}'}
\def\Dz{\bb{\DD}'_{\zz}}
\def\Ds{\bb{\DD}'_{\ss}}
\def\Dsz{\bb{\DD}'_{\so}}
\def\Dg{\bb{\DD}'_{\gg}}
\def\Ls{\bb{L}^2_{\ss}}
\def\Lg{\bb{L}^2_{\gg}}
\def\bF{{\bb{V}}}
\def\Fz{\bF_{\zz}}
\def\Fs{\bF_\ss}
\def\Fg{\bF_\gg}
\def\Pre{P}
\def\UUU{{\mathcal U}}
\def\fiapp{\phi}
\def\PU{P1}
\def\PD{P2}
\def\PT{P3}
\def\PQ{P4}
\def\PC{P5}
\def\PS{P6}
\def\Q{P6}
\def\X{Q2}
\def\Xp{Q3}
\def\Vi{V}
\def\bVi{\bb{V}}
\def\K{V}
\def\Ks{\bb{\K}_\ss}
\def\Kz{\bb{\K}_0}
\def\KM{\bb{\K}_{\, \so}}
\def\HGG{\bb{H}^\G}
\def\HG{\bb{H}^\G_{\so}}
\def\EG{{\mathfrak{P}}^{\G}}
\def\G{G}
\def\de{\delta}
\def\esp{\sigma}
\def\dd{\displaystyle}
\def\LP{\mathfrak{L}}
\def\dive{\mbox{div}}
\def\la{\langle}
\def\ra{\rangle}
\def\um{u_{\meno}}
\def\uv{\mu_{\meno}}
\def\Fp{ {\textbf F_{\piu}} }
\def\Ff{ {\textbf F} }
\def\Fm{ {\textbf F_{\meno}} }
\def\Eb{ {\textbf E} }
\def\piu{\scriptscriptstyle{+}}
\def\meno{\scriptscriptstyle{-}}
\def\omeno{\scriptscriptstyle{\ominus}}
\def\Tt{ {\mathscr T} }
\def\Xx{ {\textbf X} }
\def\Yy{ {\textbf Y} }
\def\Ee{ {\textbf E} }
\def\VP{{\mbox{\tt VP}}}
\def\CP{{\mbox{\tt CP}}}
\def\cp{$\CP(f_0, t_0)\,$}
\def\cop{$\CP(f_0)\,$}
\def\copn{$\CP_n(f_0)\,$}
\def\vp{$\VP(f_0, t_0)\,$}
\def\vop{$\VP(f_0)\,$}
\def\vopn{$\VP_n(f_0)\,$}
\def\vopdue{$\VP_2(f_0)\,$}
\def\leqs{\leqslant}
\def\geqs{\geqslant}
\def\mat{{\frak g}}
\def\tG{t_{\scriptscriptstyle{G}}}
\def\tN{t_{\scriptscriptstyle{N}}}
\def\TK{t_{\scriptscriptstyle{K}}}
\def\CK{C_{\scriptscriptstyle{K}}}
\def\CN{C_{\scriptscriptstyle{N}}}
\def\CG{C_{\scriptscriptstyle{G}}}
\def\CCG{{\mathscr{C}}_{\scriptscriptstyle{G}}}
\def\tf{{\tt f}}
\def\ti{{\tt t}}
\def\ta{{\tt a}}
\def\tc{{\tt c}}
\def\tF{{\tt R}}
\def\C{{\mathscr C}}
\def\P{{\mathscr P}}
\def\V{{\mathscr V}}
\def\TI{\tilde{I}}
\def\TJ{\tilde{J}}
\def\Lin{\mbox{Lin}}
\def\Hinfc{ H^{\infty}(\reali^d, \complessi) }
\def\Hnc{ H^{n}(\reali^d, \complessi) }
\def\Hmc{ H^{m}(\reali^d, \complessi) }
\def\Hac{ H^{a}(\reali^d, \complessi) }
\def\Dc{\DD(\reali^d, \complessi)}
\def\Dpc{\DD'(\reali^d, \complessi)}
\def\Sc{\SS(\reali^d, \complessi)}
\def\Spc{\SS'(\reali^d, \complessi)}
\def\Ldc{L^{2}(\reali^d, \complessi)}
\def\Lpc{L^{p}(\reali^d, \complessi)}
\def\Lqc{L^{q}(\reali^d, \complessi)}
\def\Lrc{L^{r}(\reali^d, \complessi)}
\def\Hinfr{ H^{\infty}(\reali^d, \reali) }
\def\Hnr{ H^{n}(\reali^d, \reali) }
\def\Hmr{ H^{m}(\reali^d, \reali) }
\def\Har{ H^{a}(\reali^d, \reali) }
\def\Dr{\DD(\reali^d, \reali)}
\def\Dpr{\DD'(\reali^d, \reali)}
\def\Sr{\SS(\reali^d, \reali)}
\def\Spr{\SS'(\reali^d, \reali)}
\def\Ldr{L^{2}(\reali^d, \reali)}
\def\Hinfk{ H^{\infty}(\reali^d, \KKK) }
\def\Hnk{ H^{n}(\reali^d, \KKK) }
\def\Hmk{ H^{m}(\reali^d, \KKK) }
\def\Hak{ H^{a}(\reali^d, \KKK) }
\def\Dk{\DD(\reali^d, \KKK)}
\def\Dpk{\DD'(\reali^d, \KKK)}
\def\Sk{\SS(\reali^d, \KKK)}
\def\Spk{\SS'(\reali^d, \KKK)}
\def\Ldk{L^{2}(\reali^d, \KKK)}
\def\Knb{K^{best}_n}
\def\sc{\cdot}
\def\k{\mbox{{\tt k}}}
\def\x{\mbox{{\tt x}}}
\def\g{ {\textbf g} }
\def\QQQ{ {\textbf Q} }
\def\AAA{ {\textbf A} }
\def\gr{\mbox{gr}}
\def\sgr{\mbox{sgr}}
\def\loc{\mbox{loc}}
\def\PZ{{\Lambda}}
\def\PZAL{\mbox{P}^{0}_\alpha}
\def\epsilona{\epsilon^{\scriptscriptstyle{<}}}
\def\epsilonb{\epsilon^{\scriptscriptstyle{>}}}
\def\lgraffa{ \mbox{\Large $\{$ } \hskip -0.2cm}
\def\rgraffa{ \mbox{\Large $\}$ } }
\def\restriction{\upharpoonright}
\def\M{{\scriptscriptstyle{M}}}
\def\m{m}
\def\Fre{Fr\'echet~}
\def\I{{\mathcal N}}
\def\ap{{\scriptscriptstyle{ap}}}
\def\fiap{\varphi_{\ap}}
\def\dfiap{{\dot \varphi}_{\ap}}
\def\DDD{ {\mathfrak D} }
\def\BBB{ {\textbf B} }
\def\EEE{ {\textbf E} }
\def\GGG{ {\textbf G} }
\def\TTT{ {\textbf T} }
\def\KKK{ {\textbf K} }
\def\HHH{ {\textbf K} }
\def\FFi{ {\bf \Phi} }
\def\GGam{ {\bf \Gamma} }
\def\sc{ {\scriptstyle{\bullet} }}
\def\a{a}
\def\ep{\epsilon}
\def\c{\kappa}
\def\parn{\par \noindent}
\def\teta{M}
\def\elle{L}
\def\ro{\rho}
\def\al{\alpha}
\def\si{\sigma}
\def\be{\beta}
\def\ga{\gamma}
\def\te{\vartheta}
\def\ch{\chi}
\def\et{\eta}
\def\complessi{{\bf C}}
\def\len{{\bf L}}
\def\reali{{\bf R}}
\def\interi{{\bf Z}}
\def\Z{{\bf Z}}
\def\naturali{{\bf N}}
\def\Sfe{ {\bf S} }
\def\To{ {\bf T} }
\def\Td{ {\To}^d }
\def\Tt{ {\To}^3 }
\def\Zd{ \interi^d }
\def\Zt{ \interi^3 }
\def\Zet{{\mathscr{Z}}}
\def\Ze{\Zet^d}
\def\T1{{\textbf To}^{1}}
\def\es{s}
\def\ee{{E}}
\def\FF{\mathcal F}
\def\FFu{ {\textbf F_{1}} }
\def\FFd{ {\textbf F_{2}} }
\def\GG{{\mathcal G} }
\def\EE{{\mathcal E}}
\def\KK{{\mathcal K}}
\def\PP{{\mathcal P}}
\def\PPP{{\mathscr P}}
\def\PN{{\mathcal P}}
\def\PPN{{\mathscr P}}
\def\QQ{{\mathcal Q}}
\def\J{J}
\def\Np{{\hat{N}}}
\def\Lp{{\hat{L}}}
\def\Jp{{\hat{J}}}
\def\Pp{{\hat{P}}}
\def\Pip{{\hat{\Pi}}}
\def\Vp{{\hat{V}}}
\def\Ep{{\hat{E}}}
\def\Gp{{\hat{G}}}
\def\Kp{{\hat{K}}}
\def\Ip{{\hat{I}}}
\def\Tp{{\hat{T}}}
\def\Mp{{\hat{M}}}
\def\La{\Lambda}
\def\Ga{\Gamma}
\def\Si{\Sigma}
\def\Upsi{\Upsilon}
\def\Gam{\Gamma}
\def\Gag{{\check{\Gamma}}}
\def\Lap{{\hat{\Lambda}}}
\def\Upsig{{\check{\Upsilon}}}
\def\Kg{{\check{K}}}
\def\ellp{{\hat{\ell}}}
\def\j{j}
\def\jp{{\hat{j}}}
\def\BB{{\mathcal B}}
\def\LL{{\mathcal L}}
\def\MM{{\mathcal U}}
\def\SS{{\mathcal S}}
\def\DD{D}
\def\Dd{{\mathcal D}}
\def\VV{{\mathcal V}}
\def\WW{{\mathcal W}}
\def\OO{{\mathcal O}}
\def\RR{{\mathcal R}}
\def\TT{{\mathcal T}}
\def\AA{{\mathcal A}}
\def\CC{{\mathcal C}}
\def\JJ{{\mathcal J}}
\def\NN{{\mathcal N}}
\def\HH{{\mathcal H}}
\def\XX{{\mathcal X}}
\def\XXX{{\mathscr X}}
\def\YY{{\mathcal Y}}
\def\ZZ{{\mathcal Z}}
\def\CC{{\mathcal C}}
\def\cir{{\scriptscriptstyle \circ}}
\def\circa{\thickapprox}
\def\vain{\rightarrow}
\def\salto{\vskip 0.2truecm \noindent}
\def\spazio{\vskip 0.5truecm \noindent}
\def\vs1{\vskip 1cm \noindent}
\def\fine{\hfill $\square$ \vskip 0.2cm \noindent}
\def\ffine{\hfill $\lozenge$ \vskip 0.2cm \noindent}
\newcommand{\rref}[1]{(\ref{#1})}
\def\beq{\begin{equation}}
\def\feq{\end{equation}}
\def\beqq{\begin{eqnarray}}
\def\feqq{\end{eqnarray}}
\def\barray{\begin{array}}
\def\farray{\end{array}}
\makeatletter \@addtoreset{equation}{section}
\renewcommand{\theequation}{\thesection.\arabic{equation}}
\makeatother
\begin{titlepage}
{~}
\vspace{-2cm}
\begin{center}
{\huge On approximate solutions of the incompressible Euler and Navier-Stokes equations}
\end{center}
\vspace{0.5truecm}
\begin{center}
{\large
Carlo Morosi$\,{}^a$, Livio Pizzocchero$\,{}^b$({\footnote{Corresponding author}})} \\
\vspace{0.5truecm} ${}^a$ Dipartimento di Matematica, Politecnico di Milano,
\\ P.za L. da Vinci 32, I-20133 Milano, Italy \\
e--mail: carlo.morosi@polimi.it \\
${}^b$ Dipartimento di Matematica, Universit\`a di Milano\\
Via C. Saldini 50, I-20133 Milano, Italy\\
and Istituto Nazionale di Fisica Nucleare, Sezione di Milano, Italy \\
e--mail: livio.pizzocchero@unimi.it
\end{center}
\begin{abstract}
We consider the incompressible Euler or Navier-Stokes (NS) equations on a torus
$\Td$, in the functional setting of the Sobolev spaces $\HM{n}(\Td)$
of divergence free, zero mean vector fields on $\Td$, for $n \in (d/2+1, +\infty)$.
We present a general theory of approximate solutions for the Euler/NS
Cauchy problem; this allows to infer a lower bound $\Tc$ on the time of existence
of the exact solution $u$ analyzing \textsl{a posteriori}
any approximate solution $\ua$, and also to construct a
function $\Rr_n$ such that $\| u(t) - \ua(t) \|_n \leqs \Rr_n(t)$ for all $t \in [0,\Tc)$.
Both $\Tc$ and $\Rr_n$ are determined solving suitable ``control inequalities'',
depending on the error of $\ua$;
the fully quantitative implementation of this scheme
depends on some previous estimates of ours on
the Euler/NS quadratic nonlinearity \cite{cok} \cite{cog}.
To keep in touch with the existing literature on the subject,
our results are compared with a setting
for approximate Euler/NS solutions proposed in \cite{Che}.
As a first application of the present framework, we consider the Galerkin approximate
solutions of the Euler/NS Cauchy problem, with a specific initial datum
considered in \cite{Nec}: in this case our methods allow, amongst else, to prove
global existence for the NS Cauchy problem when the viscosity is above
an explicitly given bound.
\end{abstract}
\vspace{0.2cm} \noindent
\hfill \par
\par \vspace{0.05truecm} \noindent \textbf{AMS 2000 Subject classifications:} 35Q30, 76D03, 76D05.
\end{titlepage}
\section{Introduction}
\label{intro}
In recent years, there has been some activity about approximate
solutions of the Euler and Navier-Stokes (NS) equations, viewed
as tools to infer accurate \textsl{a posteriori} estimates on the exact solutions.
We mention, in particular:
the works by Chernyshenko et al. \cite{Che}, Dashti and Robinson \cite{DR},
Robinson and Sadowski \cite{Rob}, and our papers \cite{uno} \cite{due} \cite{accau}.
The present work seats within the same research area; here we consider the incompressible Euler/NS equations
\beq {\partial u \over \partial t}  = - \LP(u \sc \partial u) + \nu \Delta u + f~, \label{eul} \feq
where: $u= u(x, t)$ is the divergence free velocity field;
the space variables $x = (x_s)_{s=1,...,d}$ belong to the torus
$\Td$
(and yield the derivatives $\partial_s := \partial/\partial x_s$);
$\Delta := \sum_{s=1}^d \partial_{s s}$ is the Laplacian;
$(u \sc \partial u)_r := \sum_{s=1}^d u_s \partial_s u_r$ ($r=1,...,d$);
$\LP$ is the Leray projection onto the space of divergence free vector fields;
$\nu$ is the viscosity coefficient, so that
$\nu = 0$ in the Euler case and $\nu \in (0,+\infty)$ in the NS case; $f = f(x,t)$
is the Leray projected density of external forces.
The dimension $d$ is arbitrary in the general setting of the paper,
but we put $d=3$ in a final application. \par
The functional setting that we consider
for Eq. \rref{eul} relies on the Sobolev spaces
\beq \HM{n}(\Td) \equiv \HM{n} := \{ v : \Td \vain \reali^d~|~~
 \la v \ra = 0, ~\dive v = 0,~ \sqrt{-\Delta}^{\,n} v \in \bb{L}^2(\Td)~\}~, \feq
with $\la ~ \ra$ indicating the mean over $\Td$; for any real $n$,
the above space is equipped with
the inner product $\la v | w \ra_n := \la \sqrt{-\Delta}^{\,n} v |  \sqrt{-\Delta}^{\,n} w \ra_{L^2}$
and with the corresponding norm $\| ~ \|_n$. One of the main issues in this setting
is the behavior of the bilinear map
\beq \PPP(v, w) := - \LP(v \sc w) \feq
in the above mentioned Sobolev spaces. It is well known that there are
positive constants $K_{n d} \equiv K_n$ and $G_{n d} \equiv G_n$
fulfilling the ``basic inequality''
\beq \| \PPP(v,w) \|_n \leqs K_n \| v \|_n \| w \|_{n+1} \qquad \mbox{
for $n \in (\dd{d \over 2}, + \infty)$, $v \in \HM{n}$, $w \in \HM{n+1}$}~, \label{bbasineq} \feq
and the so-called ``Kato inequality''
\beq  | \la \PPP(v,w) | w \ra_n | \leqs G_n \| v \|_n
\| w \|^{2}_{n} \quad \mbox{for $n \in (\dd{d \over 2} + 1, + \infty)$,
$v \in \HM{n}$, $w \in \HM{n+1}$}~; \label{kkatineq} \feq
fully quantitative upper and lower bounds on $K_n$ and $G_n$ were
derived in our previous works \cite{cok} \cite{cog}, for reasons
related to the present setting and described more precisely
in the sequel. \par
Independently of the problem to estimate $K_n$ and $G_n$, the above
two inequalities play a major role in the very interesting paper
\cite{Che} on approximate Euler/NS solutions and \textsl{a posteriori}
estimates on exact solutions.
To give an idea of the framework of \cite{Che} we describe a result
therein, using notations closer to our setting. \par
Consider the Euler/NS equation \rref{eul} with a specified
initial condition $u(x,0) = u_0(x)$; let $\ua : \Td \times [0,\Ta) \vain \reali^d$ be
an approximate solution of this Cauchy problem.
Given $n \in (d/2+1, +\infty)$ (and assuming suitable regularity
for $u_0, f, \ua$), let $\ua$ possess the differential error estimator
$\ep_n : [0,\Ta) \vain [0,+\infty)$, the datum error estimator $\delta_n \in [0,+\infty)$
and the growth estimators $\Dd_n, \Dd_{n+1} : [0,\Ta) \vain [0,+\infty)$; this means that,
for $t \in [0,\Ta)$,
\beq \| \big({\partial \ua \over \partial t} +
\LP(\ua \,\sc \,\partial \ua) - \nu \Delta \ua - f \big)(t)\|_n
\leqs \ep_n(t)~, \feq
\beq  \| \ua(0) - u_0 \|_n \leqs \delta_n~, \feq
\beq \| \ua(t) \|_n \leqs \Dd_n(t) ~, \quad \| \ua(t) \|_{n+1} \leqs \Dd_{n+1}(t) \feq
(with $\ua(t) := \ua(\cdot, t)$,
etc.). According to \cite{Che},
Eq. \rref{eul} with datum $u_0$ has an exact (strong, $\HM{n}$-valued) solution $u$
on a time interval $[0,\Tb) \subset [0,\Ta)$,
if $\Tb$ (with the estimators for $\ua$) fulfills the inequality
\beq \delta_n + \int_{0}^{\Tb} \!\!\! dt\, \ep_n(t)  < {1 \over G_n \Tb} ~
e^{\dd{- \!\! \!\int_{0}^{\Tb} \!\!\! dt \, \big(G_n \Dd_n(t) +  K_n \Dd_{n+1}(t)\big) }}. \label{criter} \feq
The present work aims to refine, to some extent, the approach of \cite{Che} and to apply
it to get fully quantitative estimates on the exact solution of the Euler/NS Cauchy problem on $\Td$, with
some specific initial datum. Our main
result can be described as follows: assuming
suitable regularity for $u_0, f, \ua$, and
intending $\ep_n, \Dd_n, \Dd_{n+1}$ as above, suppose
there is a function $\Rr_n \in C([0,\Tc), [0,+\infty))$,
with $\Tc \in (0,\Ta]$, fulfilling the \textsl{control inequalities}
\beq {d^{+} \Rr_n \over d t} \geqs - \nu \Rr_n
+ (G_n \Dd_n + K_n \Dd_{n+1}) \Rr_n + G_n \Rr^2_n + \ep_n
~\mbox{on $[0,\Tc)$},~~\Rr_n(0) \geqs \delta_n \label{ciq} \feq
(with $d^{+}/ d t$ the right upper Dini derivative, see Section \ref{prel}).
Then, the solution $u$ of the Euler/NS equation \rref{eul} with initial datum $u_0$
exists (in a classical sense) on the time interval $[0,\Tc)$, and its distance from the
approximate solution admits the bound
\beq \| u(t) - \ua(t) \|_n \leqs \Rr_n(t) \qquad \mbox{for $t \in [0,\Tc)$}~. \label{buap} \feq
Some features distinguishing our approach from
\cite{Che} are the following ones. \par \noindent
(i) Differently from \rref{criter}, our control inequalities \rref{ciq} depend explicitly on $\nu$
and thus could allow a more accurate analysis of the influence of viscosity
on the regularity of the Euler/NS solutions. \par \noindent
(ii) Our approach promises better lower bounds
on the time of existence of $u$. For example, for $\nu > 0$ and
under specific assumptions illustrated in the paper, the inequalities \rref{ciq} have solutions
$\Rr_n$ with $\Tc = + \infty$, implying
the global nature of the NS solution $u$;
on the contrary, if $\delta_n$ or $\ep_{n}$ are nonzero
the inequality \rref{criter} cannot have a solution with very large $\Tb$,
since the right hand side is bounded by $1/(G_n \Tb)$ and thus vanishes
for $\Tb \vain + \infty$. \par \noindent
(iii) In \cite{Che} there is not an explicit bound on
the distance between $u$ and $\ua$, such as \rref{buap} (however,
our analysis yielding \rref{buap} is greatly indebted to \cite{Che}
and, in a sense, it mainly refines and completes a chain of inequalities for
$u - \ua$ appearing therein). \par \noindent
(iv) The constants $K_n$ and $G_n$ in the inequalities
\rref{bbasineq} \rref{kkatineq} are not evaluated in \cite{Che}. On the contrary, here we have at hand
our previous results \cite{cok} \cite{cog} on these constants;
thus, in specific applications,
we can implement the control inequalities \rref{ciq} and their outcome
\rref{buap} in a fully quantitative way.
\par
As an example of our approach, in the final part of
the paper we consider the Euler/NS equations on $\Tt$
with a specific initial datum $u_0$. Independently of the approach developed
here, this datum has been already considered in an interesting paper by Behr, Ne$\check{\mbox{c}}$as
and Wu \cite{Nec}, where it is indicated as the origin of a possible blow-up
for the Euler equations.
However, in the cited work the blow-up is conjectured
on the grounds of a merely ``experimental'' analysis
of a finite number of terms in the power series $\sum_{i=0}^{+\infty}
u_i(x) t^i$ solving formally the Euler Cauchy problem. \par
In the present work, dealing with the initial datum of \cite{Nec}
both for $\nu=0$ and for $\nu > 0$, a different approach to the Cauchy
problem is developed using the familiar Galerkin approximation
(with a convenient set of Fourier modes), combined with our general
setting for approximate solutions based on the control
inequalities \rref{ciq}; in this case, the Sobolev order
is $n=3$, $\ua$ is the Galerkin solution
and we use for it the required estimators, to be substituted
in the control inequalities \rref{ciq} (with the values
for the constants $K_3$ and $G_3$ obtained
in \cite{cok} \cite{cog}). We search for
a solution $\Rr_3$ fulfilling Eqs. \rref{ciq}
as equalities (i.e., with $\leqs$ replaced
by $=$); this gives rise to an ordinary Cauchy
problem for $\Rr_3$, which is solved very easily
and reliably by numerical means.
Admittedly, our computations are preliminary:
they were performed using MATHEMATICA on a PC, with a fairly small
set of $150$ Fourier modes for the Galerkin approximation;
we plan to develop the same approach
with more powerful computational tools in a subsequent work. \par
In a few words, our results are as follows: in the
case $\nu=0$, the solution $\Rr_3$ of the control
equations \rref{ciq} exists on a finite time interval $[0, \Tc)$
(after which it blows up); so, we can grant existence for the Euler
Cauchy problem on the interval $[0, \Tc)$, where
we also have the estimate \rref{buap} on the $\HM{3}$ distance
between the exact solution $u$ and the Galerkin approximate
solution $\ua$. (Unfortunately, $\Tc$ is less than the blow-up
time suggested in \cite{Nec} for the Euler Cauchy
problem, so we cannot disprove the conjecture of the cited
paper; the situation could change using
many more Galerkin modes, which is our aim for the future).
For $0 < \nu \lesssim 8$, the situation is similar:
$\Rr_3$ blows up in a finite time $\Tc$, and
we can grant existence for the NS Cauchy problem
only up to $\Tc$. On the contrary, for $\nu \gtrsim 8$,
our approach grants global existence for the NS Cauchy
problem (and a bound of the type \rref{buap}
on the full interval $[0,+\infty)$). \par
To conclude this Introduction, let us describe
the organization of the paper. In Section \ref{prel}
we present some preliminaries: these concern mainly
the Sobolev spaces on $\Td$, in view of their applications
to the Euler/NS equations \rref{eul}. In Section \ref{caupro}
we define formally the Euler/NS Cauchy problem, the general
notion of approximate solution for this problem
and the related error estimators. In Section \ref{secmain},
that contains the main theoretical results
of the paper, we develop the general framework
yielding the control inequalities \rref{ciq},
and prove the estimate \rref{buap} on the distance
between the exact solution $u$ of the Euler/NS
Cauchy problem and an approximate solution $\ua$
(here we also give more details on the connections of the present
work with \cite{Che}).
In Section \ref{analyt} we present some analytical solutions
of the control inequalities \rref{ciq}, under specific
assumptions for their estimators and supposing,
for simplicity, that the external forcing in
\rref{eul} is zero; as anticipated, in certain cases
our analytical solutions for the control inequalities are
global, thus ensuring global existence
for the Euler/NS Cauchy problem.
In Section \ref{galegen} we describe the general Galerkin
method for \rref{eul}; in particular,
we give error and growth estimators for
the Galerkin approximate solutions, to be used
with our control inequalities \rref{ciq}.
In Section \ref{caunecas} we consider the Galerkin
method with the initial datum of \cite{Nec},
both for $\nu=0$ an for $\nu > 0$.
In Appendix \ref{appekapl} we review
some comparison lemmas of the \v{C}aplygin
type about differential inequalities; these
are employed in Section \ref{secmain} in relation
to the control inequalities. In Appendix \ref{appegal},
for completeness we report the proof
of an essentially known statement on the Galerkin approximants
for the Euler/NS Cauchy problem.
\salto
\section{Preliminaries}
\label{prel}
\noindent
\textbf{Dini derivatives.} Consider a function
\beq \fun : [0,T) \vain \reali~, \qquad t \mapsto \fun(t) \feq
(with $T \in (0,+\infty)$). The \textsl{right}, lower and upper
Dini derivatives of $\fun$ at any point $t_0 \in [0,T)$ are, respectively,
\beq {d_{+} \fun \over d t}(t_0) := \liminf_{h \vain 0^{+}} {\fun(t_0 + h) - \fun(t_0) \over h} \in [-\infty, + \infty]~;
\label{dlower} \feq
\beq {d^{+} \fun \over d t}(t_0) := \limsup_{h \vain 0^{+}} {\fun(t_0 + h) - \fun(t_0) \over h} \in [-\infty, + \infty]~.
\label{dupper} \feq
Of course,
\beq {d_{+} \fun \over d t}(t_0) \leqs {d^{+} \fun \over d t}(t_0)~; \feq
furthermore, the opposite function $-\fun : t \in [0,T) \mapsto - \fun(t)$ is such that
\beq {d_{+} (-\fun) \over d t}(t_0) = - {d^{+} \fun \over d t}(t_0)~. \feq
The \textsl{left}, lower and upper Dini derivatives $\dd{d_{-}  \over d t}$,
$\dd{d^{-}  \over d t}$ are defined similarly, with $h \vain 0^{-}$; however,
left derivatives are not used in this paper. Of course, all Dini derivatives
coincide with the usual derivative if this exists.
\salto
\textbf{Sobolev spaces of vector fields on the torus; Laplacian, Leray
projection, and so on.}
We work in any space dimension $d \in \{2,3,...\}$ (using
$r,s$ as indices in $\{1,...,d\}$).
For $a, b$ in $\complessi^d$ we put $a \sc b := \sum_{r=1}^d a_r b_r$
and write $\overline{a}$ for the complex conjugate $d$-tuple
$(\overline{a_r})$; $\complessi^d$ carries the
inner product $(a,b) \mapsto \overline{a} \sc b$
and the norm $|a| := \sqrt{\overline{a} \sc a}$. We often restrict
the previous operations to $\reali^d$.
\par \noindent
We consider the torus
$\Td$, i.e., the product of $d$ copies of $\To := \reali/(2 \pi \interi)$;
a point of $\Td$ is generically written as $x = (x_s)_{s=1,...,d}$.
In the sequel we often refer to the space
\beq \DD(\Td) \equiv \DD \feq
of the real distributions on $\Td$, and to the space
\beq \Dv(\Td) \equiv \Dv := \{ v = (v_r)_{r=1,...,d}~|~v_r \in D
~\mbox{for all $r$} \}~. \feq
Elements of $\Dv$ can be interpreted as ``generalized
functions $\Td \vain \reali^d$''; in the sequel, we call them
(distributional) vector fields on $\Td$.
$\DD$ and $\Dv$ will be equipped with their
weak topologies. For more details on distributions (and on the function spaces
mentioned in the sequel) we refer, e.g., to \cite{accau}. \par
Using distributional derivatives, we can give a meaning
to several differential operators acting
on vector fields, e.g. the Laplacian $\Delta : \Dv \vain \Dv$
and the divergence $\dive : \Dv \vain D$.
Any real distribution $w \in \DD$ has a mean
$\la w \ra := \la w, 1/(2 \pi)^d \ra \in \reali$
(the right hand side in this definition indicates the action of $w$ on the
constant test function $1/(2 \pi)^d$). Passing to a vector field
$v \in \Dv$, we can define componentwise the mean $\la v \ra \in \reali^d$.
\par
Each vector field $v \in \Dv$
has a unique (weakly convergent) Fourier series expansion
\beq v = \sum_{k \in \Zd} v_k e_k~, \quad v_k \in \complessi^d~,~~
e_k(x) := {1 \over (2 \pi)^{d/2}} \, e^{i k \sc \, x}~\mbox{for $x \in \Td$}~; \label{fs} \feq
of course, the $r$-th component of $v_k$ is $\la v_r, e_{-k} \ra$
(i.e., it equals the action of
$v_r$ on the test function $e_{-k}$).
Due to the reality of $v$, the Fourier coefficients have the property
$\overline{v_k} = v_{-k}$ for all $k \in \Zd$; one has $\la v \ra = v_0/(2 \pi)^{d/2}$. \par
In the sequel we often refer to the space of zero mean vector fields,
of the divergence free (or solenoidal) vector fields and to their
intersection; these are, respectively,
\beq \Dz := \{ v \in \Dv~|~\la v \ra = 0 \},~~
\Ds := \{ v \in \Dv~|~\dive v = 0 \},~~
\Dsz := \Ds \cap \Dz~.\feq
Elements
of $\Dz$ and $\Ds$ are characterized, respectively, by the
conditions $v_0 =0$ and $k \sc v_k=0$ for all $k$;
handling with the Fourier components of
vector fields in $\Dz$, it is convenient to put
\beq \Zd_0 :=  \Zd \setminus \{0\}~. \feq
Of course, for each $v \in \Dv$ one has $(\Delta v)_k =-|k|^2 v_k$;
this suggests to define, for any $n \in \reali$,
\beq \sqrt{-\Delta}^n : \Dz \vain \Dz,~~
v \mapsto \sqrt{-\Delta}^n~\mbox{such that $(\sqrt{-\Delta}^n v)_k := |k|^n v_k$
for $k \in \Zd_0$}. \feq
We denote with $\bb{L}^2(\Td)$ the space of square integrable
vector fields $v : \Td \vain \reali^d$; this is a real Hilbert space with the inner product
$\la v | w \ra_{L^2} := \int_{\Td} d x \, v(x) \sc w(x)$ $=\sum_{k \in \Zd} \overline{v_{k}} \sc w_k$
and the norm $\| v \|_{L^2} := \sqrt{\la v | v \ra_{L^2}}$.
For any $n \in \reali$, let us consider the Sobolev space
\beq \Hz{n}(\Td) \equiv \Hz{n} := \{ v \in \Dz~|~
\sqrt{-\Delta}^{\,n} v \in \bb{L}^2(\Td)~\} = \feq
$$ = \{ v \in \Dv~|~
v_0 = 0, \sum_{k \in \Zd_0} |k|^{2 n} |v_k|^2 < + \infty~\}~; $$
this is a real Hilbert space with the inner product and the norm
\beq \la v | w \ra_n := \la \sqrt{-\Delta}^{\,n} v |  \sqrt{-\Delta}^{\,n} w \ra_{L^2} =
\sum_{k \in \Zd_0} |k|^{2 n} \overline{v_k} \sc w_k~, \feq
\beq \| v \|_n := \sqrt{\la v | v \ra_n} = \sqrt{\sum_{k \in \Zd_0} |k|^{2 n} |v_k|^2}~. \feq
Clearly, $n \leqs n'$ implies $\Hz{n'} \subset \Hz{n}$ and
$\|~\|_{n} \leqs \|~\|_{n'}$.
In this paper, we mainly fix the attention on the divergence free Sobolev space
\beq \HM{n}(\Td) \equiv \HM{n} := \Ds \cap \Hz{n} = \{ v \in \Hz{n}~|~k \sc v_k = 0~\} \feq
($n \in \reali$); this is a closed subspace of $\Hz{n}$, and thus a real Hilbert space with
the restriction of $\la~|~\ra_n$. \par
For $m \in \{0,1,2,...,+\infty\}$ let us
consider the space $\mathbb{C}^m(\Td, \reali^d) \equiv \mathbb{C}^m$
of vector fields $v : \Td \vain \reali^d$ of class $C^m$ in the ordinary sense;
if $m < + \infty$, equip this with the norm
\beq \| v \|_{C^m} := \max_{\tiny{\barray{cc} \ell =0,...,m, \\ r,s_1,...,s_\ell=1,...,d \farray}}
\max_{x \in \Td} |\partial_{s_1...s_\ell} u_r(x)|~. \feq
There is a well known Sobolev imbedding
\beq \Hz{n} \subset \mathbb{C}^m,~\| ~\|_{C^m} \leqs S_{m n} \|~\|_n
\quad \mbox{for $m \in \{0,1,2,...\}$, $n \in \reali$, $n  > \dd{d \over 2} + m$}~, \label{imb} \feq
involving suitable positive constants $S_{m n d} \equiv S_{m n}$. \par
For arbitrary $n \in \reali$, we have
\beq \Delta \Hz{n+2} = \Hz{n}~, \qquad \Delta \HM{n+2} = \HM{n}~, \label{defa}\feq
and $\Delta$ is continuous between the above spaces; furthermore,
\beq \la \Delta v | v \ra_n = - \| v \|^2_{n+1}
\leqs - \| v \|^2_n \qquad \mbox{for each $v \in \Hz{n+2}$}
\label{inelap} \feq
(all the above statements are made evident by the Fourier representations). \par
By definition, the \textsl{Leray projection} is the map
\beq \LP : \Dv \vain \Ds~, \qquad v \mapsto \LP v~\mbox{such that
$(\LP v)_k = \LP_k v_k$ for $k \in \Zd$}~; \feq
here $\LP_k$ is the orthogonal projection of $\complessi^d$ onto $k^\perp = \{ a \in
\complessi^d~|~k \sc a = 0 \}$ (given explicitly by $\LP_k c = c
- (k \sc c) k/|k|^2$ and $\LP_0 c = c$, for $k \in \Zd_0$ and
$c \in \complessi^d$).
For each real $n$, we have a continuous linear map
\beq \LP \restriction \Hz{n} : \Hz{n} \vain \HM{n} \feq
(which, in fact, is the orthogonal projection of
$(\Hz{n}, \la ~|~\ra_n)$ onto the subspace
$\HM{n}$).
\salto
\textbf{The fundamental bilinear map in the Euler/NS equations.}
In this subsection we assume $n \in (d/2,+ \infty)$. With this condition, we have
a continuous bilinear map
$\HM{n} \times \HM{n+1} \vain \Hz{n}$, $(v, w) \mapsto v \sc \partial w$,
where  $v \sc \partial w$ is the vector field on $\Td$ with components $(v \sc \partial w)_r :=
\sum_{s = 1}^d v_s \partial_s w_r$. By composition with the Leray projection $\LP$,
we obtain a continuous bilinear map
\beq \PPP : \HM{n} \times \HM{n+1} \vain \HM{n}~, \qquad (v, w) \mapsto \PPP(v,w) := - \LP(v \sc \partial w)~. \feq
This appears in the Euler/NS equations, and is referred in the sequel
as the fundamental bilinear map for such equations.
As is known, for all $v \in \HM{n}$ and $w \in \HM{n+1}$ of Fourier components $v_k$ and $w_k$,
$\PPP(v, w)$ has Fourier components $\PPP(v,w)_0 = 0$ and
\beq \PPP(v, w)_k = - {i \over (2 \pi)^{d/2}}
\sum_{h \in \Zd_0} v_h \sc (k-h) \, \LP_k w_{k-h} \label{repp} \feq
for all $k \in \Zd_0$, where $\LP_k$ is the already mentioned
projection of $\complessi^d$ onto $k^\perp$.
\par
The continuity of $\PPP$ is equivalent to the existence of a constant $K_{n d} \equiv K_n \in (0,+\infty)$ such that
\beq \| \PPP(v,w) \|_n \leqs K_n \| v \|_n \| w \|_{n+1} \qquad \mbox{
for $v \in \HM{n}$, $w \in \HM{n+1}$}~; \label{basineq} \feq
we refer to this as the \textsl{basic inequality} about $\PPP$.
With the stronger assumption $n \in (d/2 + 1, +\infty)$, it is known
that there is a constant $G_{n d} \equiv G_n \in (0,+\infty)$ such that
\beq  | \la \PPP(v,w) | w \ra_n | \leqs G_n \| v \|_n
\| w \|^{2}_{n} \quad \mbox{for $v \in \HM{n}$, $w \in \HM{n+1}$}~; \label{katineq} \feq
we call this the \textsl{Kato inequality}, since it originates from
Kato's seminal paper \cite{Kato} (for completeness, we mention that
$\la \PPP(v,w) | w \ra_n = - \la v \sc \partial w | w \ra_n$ for $v, w$ as above).
In our previous works \cite{cok} \cite{cog}, we
derived upper and lower bounds for the sharp constants in the above inequalities;
throughout this paper, $K_n$ and $G_n$ are any two constants fulfilling
Eqs. \rref{basineq} \rref{katineq}.
\section{The Cauchy problem for the Euler/NS equations: exact and
approximate solutions}
\label{caupro}
From here to the end of Section \ref{galegen}
we fix any space dimension $d \in \{2,3,...\}$, we consider the Sobolev spaces
of vector fields on $\Td$, and we choose a real number $n$
such that
\beq n \in ({d \over 2} + 1, + \infty)~. \label{propen} \feq
\textbf{Euler and NS equations: the Cauchy problem.}
Let us choose a ``viscosity coefficient''
\beq \nu \in [0,+\infty)~, \label{vis}\feq
a ``forcing''
\beq f \in C([0,+\infty), \HM{n}) \label{for}\feq
and an initial datum
\beq u_0 \in \HM{n+2}~. \label{in}\feq
\begin{prop}
\textbf{Definition.} The Cauchy problem for the (incompressible) fluid
with viscosity $\nu$, initial datum $u_0$ and forcing $f$ is the following:
\beq \mbox{\textsl{Find}}~
u \in C([0, T), \HM{n+2}) \cap C^1([0,T), \HM{n}) \quad \mbox{\textsl{such that}} \label{cau} \feq
$$ {d u \over d t} = \nu \Delta u + \PPP(u,u) + f~, \qquad u(0) = u_0 $$
(with $T \in (0, + \infty]$, depending on $u$). As usually, we
speak of the ``Euler Cauchy problem'' if $\nu=0$, and
of the ``NS Cauchy problem'' if $\nu >0$.
\ffine
\end{prop}
\begin{rema} \textbf{Remark.} (i) The map $u \mapsto \PPP(u,u)$ sends continuously $\HM{n+1}$
into $\HM{n}$, while $\Delta$ sends continuously $\HM{n+2}$ into $\HM{n}$;
this explains the appearing of $\HM{n+2}$ in the previous definition,
at least in the NS case. In the Euler case $\nu=0$, $\Delta$
is absent from \rref{cau}; so, in the previous definition and in
the subsequent theoretical developments one could systematically
replace $\HM{n+2}$ with $\HM{n+1}$; this is not done just
to avoid tedious distinctions between the Euler and the NS case.
\parn
(ii) Consider any function $u \in C([0,T),\Hz{n+2})$. By the Sobolev imbedding
\rref{imb}, one has
$$ u \in C([0,T),\mathbb{C}^m) \quad
\mbox{for any $m \in \{0,1,2,...\}$ such that $n+2 > \dd{d \over 2} + m$} $$
\beq \mbox{(in particular, for $m=0,1,2,3$)}~. \label{obprp} \feq
\ffine
\end{rema}
The following results are well known, and reported for completeness.
\begin{prop}
\textbf{Proposition.} For any $\nu, f$ and $u_0$ as above, (i)-(iii) hold.
\par \noindent
(i) The Cauchy problem \rref{cau} has a unique maximal (i.e.,
nonextendable) solution, hereafter denote with $u$,
with a suitable domain $[0,T)$ ($0 < T \leqs +\infty$).
All the other solutions of \rref{cau} are restrictions
of $u$. \parn
(ii) $u$ has the property \rref{obprp}.
Furthermore, if $\nu > 0$ and the forcing $(x,t) \mapsto f(t)(x) \equiv f(x,t)$ is $C^\infty$
from $\Td \times (0,T)$ to $\reali^d$, the function $(x, t) \mapsto u(t)(x) \equiv u(x,t)$
is $C^\infty$ as well from $\Td \times (0, T)$ to $\reali^d$. \parn
(iii) For any $t \in [0,T)$ consider the vorticity matrix
$\omega(t)$, of elements
$\omega_{r s}(t) := \partial_r u_s(t) - \partial_s u_r(t)$
($r,s=1,...,d$);
let $\| \omega(t) \|_{C^0} := \max_{r,s=1,...,d} \max_{\Td}|\omega_{r s}(t)|$.
If $T < +\infty$, one has
\beq \limsup_{t \vain T^{-}} \| \omega(t) \|_{C^0} = + \infty~; \label{vor} \feq
this implies
\beq \limsup_{t \vain T^{-}} \| u(t) \|_{n} = + \infty~. \label{uninf} \feq
\end{prop}
\textbf{Proof.} (i) See \cite{Kat2}. \parn
(ii) The property \rref{obprp} holds because $u \in C([0,T), \HM{n+2})$.
For $\nu > 0$ and $f$ of class $C^\infty$ on $\Td \times (0,T)$,
the same smoothness property is granted for $u$
by Theorem 6.1 of \cite{Kat3} (with the domain $D$ considered
therein replaced by $\Td$). \parn
(iii) Let $T < + \infty$.
Eq. \rref{vor} is the celebrated Beale-Kato-Majda blow-up criterion \cite{BKM} (see
also \cite{Koz}).
To prove \rref{uninf} we note that, for each $t \in [0,T)$,
\beq \| \omega(t) \|_{C^0} \leqs \mbox{const.} \| u(t) \|_{C^1} \leqs
\mbox{const.} \| u(t) \|_n~; \label{notobs} \feq
here the first inequality follows from the definition of $\omega(t)$
in terms of the derivatives $\partial_r u_s(t)$, and the second one
from the Sobolev imbedding \rref{imb} (with $m=1$).
The relations \rref{vor} and \rref{notobs} immediately give
Eq. \rref{uninf}. \fine
\textbf{Approximate solutions.} Our treatment uses systematically the present terminology.
\begin{prop}
\label{defapp}
\textbf{Definition.} An approximate solution of the problem \rref{cau} is
any map $\ua \in C([0, \Ta), \HM{n+2}) \cap C^1([0,\Ta), \HM{n})$ (with $\Ta \in (0,+\infty]$).
Given such a function, we stipulate (i) (ii). \par\noindent
(i) The differential
error of $\ua$ is
\beq e(\ua) := {d \ua \over d t} - \nu \Delta \ua - \PPP(\ua,\ua) - f~
\in C([0,\Ta), \HM{n})~;  \feq
the datum error is
\beq \ua(0) - u_0 \in \HM{n+2}~. \feq
(ii) Let $m \in \reali, m \leqs n$. A \textsl{differential error estimator of order $m$} for $\ua$ is a function
\beq \ep_m \in C([0,\Ta), [0,+\infty))
\quad \mbox{such that}\quad \| e(\ua)(t) \|_m \leqs \ep_m(t)~\mbox{for $t \in [0,\Ta)$}~. \feq
Let $m \in \reali$, $m \leqs n + 2$.
A \textsl{datum error estimator of order $m$} for $\ua$ is a real number
\beq \delta_m \in [0,+\infty) \quad \mbox{such that} \quad \| \ua(0) - u_0 \|_m \leqs \delta_m~; \feq
a \textsl{growth estimator of order $m$} for $\ua$ is a function
\beq \Dd_m \in C([0,\Ta), [0,+\infty))
\quad \mbox{such that} \quad \| \ua(t) \|_m \leqs \Dd_m(t)~\mbox{for $t \in [0,\Ta)$}~. \label{din} \feq
In particular the function $\ep_m(t) := \| e(\ua)(t) \|_m$,
the number $\delta_m := \| \ua(0) - u_0 \|_m$ and
the function $\Dd_m(t) := \| \ua(t) \|_m$ will be called the
\emph{tautological estimators} of order $m$ for the differential error,
the datum error and the growth of $\ua$.
\end{prop}
\section{Main theorems about approximate solutions}
\label{secmain}
\textbf{Assumptions and notations.}
Throughout this section we fix
a viscosity coefficient, a forcing and an initial datum as in Eqs. \rref{vis}--\rref{in}
(recalling the condition \rref{propen} $n > d/2+1$).
We consider for the Cauchy problem \rref{cau} an approximate solution
\beq \ua \in C([0, \Ta), \HM{n+2}) \cap C^1([0,\Ta), \HM{n})~; \feq
in the sequel $\ep_n \in C([0,\Ta), [0,+\infty)$ and $\delta_n \in [0,+\infty)$
are differential error and datum error estimators of order $n$ for $\ua$, while
$\Dd_n, \Dd_{n+1} \in C([0,\Ta), [0,+\infty))$ are growth estimators of orders $n, n+1$
(see Definition \ref{defapp}).
Finally, we denote with
\beq u \in C([0, \Tm), \HM{n+2}) \cap C^1([0,\Tm), \HM{n})  \feq
the maximal solution of \rref{cau} (typically unknown, as well as $\Tm$).
\salto
\textbf{Some lemmas.} For the sake of brevity, we put
\beq \Tw := \min(\Ta, \Tm)~, \qquad w := u - \ua \in C([0, \Tw), \HM{n+2}) \cap C^1([0,\Tw), \HM{n})~.
\label{tw} \feq
Furthermore, we introduce the function
\beq \Ww_n : [0,\Tw) \vain [0,+\infty)~, \qquad \Ww_n(t) := \| w(t) \|_n~; \label{defwn} \feq
this is clearly continuous, and $C^1$ in a neighborhood of any instant $t_0 \in [0,\Tw)$
such that $w(t_0) \neq 0$. In the sequel we often consider the right, upper Dini derivative $d^{+} \Ww_n/d t$
(see Eq. \rref{dupper}), that is just the ordinary derivative at any $t_0$ with $w(t_0) \neq 0$.
\par
The forthcoming two lemmas review and partially refine some results of \cite{Che}
(after adaptation to our slightly different setting; for example, Dini
derivatives are not even mentioned in \cite{Che}).
\begin{prop}
\textbf{Lemma.} One has
\beq {d w \over d t} = \nu \Delta w + \PPP(\ua, w) + \PPP(w, \ua) + \PPP(w,w) - e(\ua)
\quad \mbox{on $[0,\Tw)$}~, \label{ethen} \feq
\beq w(0) = u_0 - \ua(0)~. \label{ew0} \feq
\end{prop}
\textbf{Proof.} By definition of the differential error $e(\ua)$, we have
\beq {d \ua \over d t} = \nu \Delta \ua + \PPP(\ua, \ua) + f + e(\ua)~; \label{wehaveuno} \feq
of course, $d u/ d t = \nu \Delta u + \PPP(u, u) + f$ whence, writing
$u = \ua+ w$,
$$ {d \ua \over d t} + {d w \over d t} = \nu \Delta (\ua + w)  + \PPP(\ua + w, \ua + w) + f  $$
\beq = \nu \Delta \ua + \nu \Delta w  + \PPP(\ua, \ua) + \PPP(\ua, w) + \PPP(w, \ua) + \PPP(w,w) + f~.
\label{eqlast} \feq
Subtracting Eq. \rref{wehaveuno} from \rref{eqlast} we obtain Eq. \rref{ethen}; Eq. \rref{ew0} is obvious. \fine
\begin{prop}
\label{lemwn}
\textbf{Lemma.} Consider the above mentioned estimators $\ep_n, \delta_n, \Dd_n, \Dd_{n+1}$,
and the function $\Ww_n$ defined in Eq. \rref{defwn}. Then
\beq {d^{+} \Ww_n \over d t} \leqs - \nu \, \Ww_n + (G_n \Dd_n + K_n \Dd_{n+1}) \Ww_n + G_n \Ww^2_n + \ep_n
~\mbox{everywhere on $[0,\Tw)$}, \label{ver1} \feq
\beq \Ww_n(0) \leqs \delta_n ~. \label{ver2} \feq
\end{prop}
\textbf{Proof.} We proceed in three steps. \par
\textsl{Step 1. Verification of Eq. \rref{ver1} at a time $t_0$ such that $w(t_0) \neq 0$.}
In this case, $w(t) \neq 0$ for all $t$ in some interval $I \ni t_0$; $\Ww_n$ is nonzero and
$C^1$ on $I$. In the same interval, we have:
\par \noindent
\vbox{
\beq {d^{+} \Ww_n \over d t} = {d \Ww_n \over d t} = {1 \over 2 \Ww_n} {d \Ww^2 _n\over d t} = {1 \over 2 \Ww_n}
{d \over d t} \la w | w \ra_n = {1 \over \Ww_n} \la {d w \over d t} | w \ra_n \label{wehave} \feq
$$ \underset{\rref{ethen}}{=} {1 \over \Ww_n} \Big(\nu \la \Delta w | w \ra_n + \la \PPP(\ua, w)|w\ra_n +
\la \PPP(w, \ua) | w \ra_n + \la \PPP(w,w) | w \ra_n - \la e(\ua) | w \ra_n \Big)~. $$
}
On the
other hand, we have the following inequalities:
\beq \la \Delta w | w \ra_n \underset{\rref{inelap}}{\leqs} - \| w \|^2_n = - \Ww^2_n~; \label{fol1} \feq
\beq \la \PPP(\ua, w)|w\ra_n \leqs | \la \PPP(\ua, w)|w\ra_n |
\underset{\rref{katineq}}{\leqs} G_n \| \ua \|_n
\| w \|^2_n \leqs G_n \Dd_n \Ww^2_n~; \feq
\vbox{
\beq \la \PPP(w,\ua)|w\ra_n \leqs | \la \PPP(w,\ua)|w\ra_n | \leqs \| \PPP(w, \ua) \|_n \|w \|_n \feq
$$ \underset{\rref{basineq}}{\leqs} K_n  \| \ua \|_{n+1} \| w \|^2_n \leqs K_n \Dd_{n+1} \Ww^2_n~; $$
}
\beq \la \PPP(w,w)|w\ra_n \leqs | \la \PPP(w,w)|w\ra_n | \underset{\rref{katineq}}{\leqs} G_n \| w \|^3_n
= G_n \Ww^3_n~; \feq
\beq - \la e(\ua) | w \ra_n \leqs |\la e(\ua) | w \ra_n| \leqs \| e(\ua) \|_n \| w \|_n \leqs \ep_n \Ww_n~;
\label{fol4} \feq
inserting the relations \rref{fol1}--\rref{fol4} into Eq. \rref{wehave}, we obtain the inequality
\rref{ver1} at all times $t \in I$ and, in particular, at time $t_0$. \par
\textsl{Step 2. Verification of Eq. \rref{ver1} at a time $t_0$ such that $w(t_0) = 0$.}
At this instant we have $\Ww_n(t_0) = 0$, and the relation \rref{ver1} to be proved becomes
\beq  {d^{+} \Ww_n \over d t}(t_0) \leqs \ep_n(t_0)~. \label{toprove} \feq
To go on we recall that, given a $C^1$ function $z: [0,T_z) \vain \Eb$ with values in a Banach space
$\Eb$ equipped with the norm $\|~\|$, we have  $d^{+} \| z \| / d t \leqs \| d z/ d t \|$
at all times, including instants when $z$ vanishes (see, e.g., \cite{Pet} and references
therein); if this result is applied
to the function $w$ with norm $\| w \|_n = \Ww_n$, we get
\beq {d^{+} \Ww_n \over d t}(t_0)  \leqs \| {d w \over d t}(t_0) \|_n~. \label{now1} \feq
Now, we evaluate $(d w/ d t)(t_0)$ via Eq. \rref{ethen}, taking into account that $w(t_0) = 0$;
this gives
\beq \| {d w \over d t}(t_0) \|_n = \| e(\ua)(t_0) \|_n \leqs \ep_n(t_0)~, \label{now2} \feq
and Eqs. \rref{now1} \rref{now2} yield the thesis \rref{toprove}. \par
\textsl{Step 3. Verification of Eq. \rref{ver2}.} In fact, the definitions of $\Ww_n, w$ and $\delta_n$ give
\par \noindent
\vbox{
\beq \Ww_n(0) = \| u(0) - \ua(0) \| \leqs \delta_n~. \feq
Now, the proof is concluded.
\fine
}
To go on, we need a comparison lemma of the \v{C}aplygin type. This is a
variant of some known results \cite{Las} \cite{Mitr}: for more details,
see Appendix \ref{appekapl}.
\begin{prop}
\label{lemtesto}
\textbf{Lemma.} Let $\Tl \in (0,+\infty]$;
consider a function $\ef \in C(\reali \times [0,\Tl),\reali)$, $(s,t)
\mapsto \ef(s,t)$ possessing partial derivative $\partial \ef/\partial s \in C(\reali \times
[0,\Tl), \reali)$;
furthermore, let $s_0 \in \reali$.
Suppose that $\Ww, \Rr \in C([0,\Tl), \reali)$ are functions such that
\beq {d_{+} \Ww \over d t}(t) \leqs \ef(\Ww(t),t)\quad\mbox{for all $t \in [0,\Tl)$},~ \Ww(0) \leqs s_0~, \feq
\beq {d^{+} \Rr \over d t}(t) \geqs \ef(\Rr(t),t)\quad\mbox{for all $t \in [0,\Tl)$},~ \Rr(0) \geqs s_0~. \feq
Then
\beq \Ww(t) \leqs \Rr(t) \qquad \mbox{for all $t \in [0,\Tl)$}~. \label{eqtesto} \feq
\ffine
\end{prop}
Hereafter, interesting conclusions will arise from
the combination of the above comparison result
with Lemma \ref{lemwn}.
\salto
\textbf{The control inequalities and the main theorem.} We begin with a definition,
which is followed by the main result of the section.
\begin{prop}
\textbf{Definition.} Consider a function $\Rr_n \in C([0,\Tc), [0,+\infty))$,
with $\Tc \in (0,\Ta]$. This function is said to fulfill the
\textsl{control inequalities} (with respect to the estimators $\ep_n, \delta_n, \Dd_n,
\Dd_{n+1}$) if
\beq {d^{+} \Rr_n \over d t} \geqs - \nu \Rr_n
+ (G_n \Dd_n + K_n \Dd_{n+1}) \Rr_n + G_n \Rr^2_n + \ep_n
~\mbox{everywhere on $[0,\Tc)$}, \label{cont1} \feq
\beq \Rr_n(0) \geqs \delta_n~. \label{cont2} \feq
\end{prop}
\begin{prop}
\label{main}
\textbf{Proposition.} Suppose there is a function $\Rr_n \in C([0,\Tc), [0,+\infty))$ fulfilling
the control inequalities, and consider the maximal solution $u$ of the
Euler/NS Cauchy problem \rref{cau}. Then, the existence time $T$ of $u$ is such that
\beq T \geqs \Tc~; \label{tta} \feq
furthermore,
\beq \| u(t) - \ua(t) \|_n \leqs \Rr_n(t) \qquad \mbox{for $t \in [0,\Tc)$}~. \label{furth} \feq
\end{prop}
\textbf{Proof.} Let us employ the following notation, already used in the previous subsections:
$\Tw := \min(\Ta, T)$, $w := u - \ua$
and $\Ww_n : t \mapsto  \| w(t) \|_n$ (see Eqs. \rref{tw} \rref{defwn};
$w$ and $\Ww_n$ have domain $[0, \Tw)$). We further define
$\Tr := \min(\Tc, T)$ (and note that $\Tc \leqs \Ta$ implies $\Tr \leqs \Tw$).
For the moment we have not yet proved \rref{tta}, so we do not know whether $\Tr$ equals $\Tc$, or not.
We proceed in three steps. \par
\textsl{Step 1. One has}
\beq \| u(t) - \ua(t) \|_n \leqs \Rr_n(t) \qquad \mbox{for $t \in [0,\Tr)$}~. \label{urw} \feq
Consider the functions $\Ww_n$ and $\Rr_n$ on $[0,\Tr)$. Due to Lemma \ref{lemwn}
and to the control inequalities, these fulfill the relations
\par\noindent
\vbox{
$$ {d_{+} \Ww_n \over d t}(t) \leqs {d^{+} \Ww_n \over d t}(t)
\leqs \ef(\Ww_n(t),t)~\mbox{for all $t \in [0,\Tr)$},~~ \Ww_n(0) \leqs \delta_n~; $$
$$ {d^{+} \Rr_n \over d t}(t) \geqs \ef(\Rr_n(t),t)~\mbox{for all $t \in [0,\Tr)$},~~ \Rr_n(0) \geqs \delta_n~, $$
where $\ef : [0,\Tr) \times \reali \vain \reali$, $\ef(s,t) :=
- \nu s + (G_n \Dd_n(t) + K_n \Dd_{n+1}(t)) s + G_n s^2 + \ep_n(t)$.
}
Therefore, by Lemma \ref{lemtesto} we have $\Ww_n(t) \leqs \Rr_n(t)$ for
$t \in [0,\Tr)$; this is just the relation \rref{urw}. \par
\textsl{Step 2. One has the relation \rref{tta} $T \geqs \Tc$}.
To prove this, we assume $T < \Tc$ and try to infer a contradiction.
To this purpose we note that, according to \rref{uninf},
$$ \limsup_{t \vain T^{-}} \| u(t) \|_n = + \infty~. $$
On the other hand, recalling Eq. \rref{urw} we can write
$\| u(t) \|_n$ $\leqs \| \ua(t) \|_n$ $ + \| u(t) - \ua(t) \|_n$ $\leqs \| \ua(t)\|_n + \Rr_n(t)$
for $t \in [0,T)$; from here and from the continuity of $\| \ua(~) \|_n$, $\Rr_n$
on $[0,\Tc) \supset [0,T]$
we get
\beq \sup_{t \in [0,T)} \| u(t) \|_n \leqs
\max_{t \in [0,T]}  (\,\| \ua(t) \|_n + \Rr_n(t)\,) < + \infty~. \label{wehave2} \feq
The results \rref{uninf} and \rref{wehave2} are contradictory. \par
\textsl{Step 3. Conclusion of the proof}. Eq. \rref{tta} proved in Step 2 implies
$\Tr = \Tc$; now, the inequality \rref{urw} of Step 1 coincides with the thesis \rref{furth}. \fine
Of course, the control inequalities \rref{cont1} \rref{cont2} are fulfilled by a function
$\Rr_n \in C^1([0,\Tc),[0,+\infty))$ such that $\Tc \in (0, \Ta]$, and
\par
\vbox{
\beq {d \Rr_n \over d t} = - \nu \Rr_n
+ (G_n \Dd_n + K_n \Dd_{n+1}) \Rr_n + G_n \Rr^2_n + \ep_n
\quad\mbox{everywhere on $[0,\Tc)$}, \label{econt1} \feq
\beq \Rr_n(0) = \delta_n~. \label{econt2} \feq
}
\noindent
\begin{prop}
\textbf{Definition.}
Eqs. \rref{econt1} \rref{econt2} are referred in the sequel as
the the \textsl{control equations}, or the \textsl{control Cauchy problem} for $\Rr_n$.
\ffine
\end{prop}
\section{Analytic solutions of the control inequalities (with no
external forcing)}
\label{analyt}
Throughout this section we consider again the Euler/NS Cauchy problem
\rref{cau}, for a given viscosity  $\nu \in [0,+\infty)$, external
forcing $f$ and initial datum $u_0$ as in
\rref{vis}--\rref{in}.
After choosing an approximate solution $\ua$ for \rref{cau} and determining the
corresponding estimators, one is faced with the problem of solving
the control inequalities \rref{cont1} \rref{cont2}, or the control
equations \rref{econt1} \rref{econt2}, for the unknown
real function $t \mapsto \Rr_n(t)$.
In many applications, such as the one of the next section on the Galerkin
approximate solutions, a numerical treatment of the control equations
is recommended (and the result is generally reliable:
\rref{econt1}  \rref{econt2} is a typically nonstiff
Cauchy problem).
However, an analytic approach to the control equations/inequalities has
its own interest, both for theoretical reasons and
for building simple user-ready criteria. \par
In this section we propose an analytic approach for special
forms of the approximate solution and/or its estimators.
For simplicity, throughout the section we assume zero external forcing:
\beq f(t) := 0 \qquad \mbox{for $t \in [0,+\infty)$}~; \feq
however, many results presented in the section could be extended
to the case of nonzero $f$, with suitable assumptions on this function. \par
We begin by considering the approximate solution
$\ua := 0$; this choice seems to be very trivial but, in spite of this,
it can be used to obtain nontrivial estimates on the time of existence
of the exact solution of the Euler/NS Cauchy problem.
In the NS case $\nu > 0$, these estimates include a condition for global
existence under a fully quantitative norm
bound on the initial datum. \par
The second case considered in the section
is much more general: the approximate solution is unspecified,
even though a certain form is assumed for its estimators. \par
\salto
\textbf{Results from the zero approximate solution.}
Let us choose for \rref{cau} the approximate solution
\beq \ua : [0,+\infty) \vain \HM{n+2}~, \qquad \ua(t) := 0 \quad \mbox{for all $t$}~. \feq
\begin{prop}
\label{lemzero}
\textbf{Lemma.} (i) The differential and datum errors of the zero approximate solution are
\beq e(\ua)(t) = 0 \qquad \mbox{for all $t \in [0,+\infty)$}~, \qquad \ua(0) - u_0 = - u_0~; \feq
consequently, this approximate solution has the differential and datum estimators
\beq \ep_n(t) := 0~, \qquad \delta_n := \| u_0 \|_n~, \label{estia} \feq
and the growth estimators
\beq \Dd_n(t) := \Dd_{n+1}(t) := 0 \qquad \mbox{for $t \in [0,+\infty)$}~. \label{estib} \feq
The control equations \rref{econt1} \rref{econt2} with the above estimators take the form
\beq {d \Rr_n \over d t} = - \nu \Rr_n + G_{n} \Rr_n^2~, \qquad \Rr_n(0) = \| u_0 \|_n~.
\label{estic} \feq
(ii) The Cauchy problem \rref{estic} has the solution $\Rr_n \in C^1([0,\Tc), [0,+\infty))$,
determined as follows:
\beq \Tc := \left\{ \barray{ll} + \infty & \mbox{if $\nu >0$, $\| u_0 \|_n \leqs {\nu/G_n}$}~, \\
- \dd{1 \over \nu} \log\left(1 - \dd{\nu \over G_n \| u_0 \|_n} \right) & \mbox{if $\nu > 0$,
$\| u_0 \|_n > {\nu/G_n}$,} \\
\dd{1 \over G_n \| u_0 \|_n} & \mbox{if $\nu=0$} \farray \right. \label{ta} \feq
(intending $1/(G_n \| u_0 \|_n) := + \infty$ if $u_0 = 0$);
\beq \Rr_n(t) := {\| u_0 \|_n e^{-\nu t} \over 1 - G_n \| u_0 \|_n e_{\nu}(t)} \qquad \mbox{for
$t \in [0,\Tc)$}~, \label{ern} \feq
\beq e_{\nu}(t) := \left\{ \barray{ll} \dd{1 - e^{-\nu t} \over \nu} & \mbox{if $\nu > 0$}, \\
t & \mbox{if $\nu = 0$} \farray \right. \label{enu} \feq
(note that $t = \lim_{\nu \vain 0^{+}} \dd{1 - e^{-\nu t} \over \nu}$).
\end{prop}
\textbf{Proof.} (i) Obvious. \par \noindent
(ii) The Cauchy problem \rref{estic} is solved by the quadrature formula \beq
\int_{\dd{\|u_0\|_n}}^{\dd{\Rr_n(t)}} {d r \over G_n r^2 - \nu r}
= t~; \feq the integral in the left hand side equals $(1/\nu)
\log \dd{G_n - \nu/\Rr_n(t) \over G_n - \nu/\|u_0 \|_n}$ if $\nu >
0$, and $(1/G_n)(1/\| u_0 \|_n - 1/\Rr_n(t))$ if $\nu=0$; some
elementary manipulations yield for the solution $\Rr_n$ the
expression \rref{ta}-\rref{enu}. \fine
From the previous lemma and from the main theorem
on approximate solutions (Proposition \ref{main}), here applied with $\ua = 0$
and $\Rr_n$ as in \rref{ern}, we obtain
the following result.
\begin{prop}
\label{prozero}
\textbf{Proposition.}
Consider the Cauchy problem \rref{cau} for the Euler/NS equations
with zero external forcing, and any datum $u_0 \in \HM{n+2}$;
let $u \in C([0, T), \HM{n+2}) \cap C^1([0,T), \HM{n})$ be
the maximal solution.
Define $\Tc$ and $e_{\nu}$ as in Eqs. \rref{ta} \rref{enu}; then
\beq T \geqs \Tc~, \qquad \| u(t) \|_n \leqs
{ \| u_0 \|_n e^{-\nu t} \over 1 - G_n \| u_0 \|_n e_{\nu}(t)} \qquad \mbox{for $t \in [0,\Tc)$}~.
\label{cobg} \feq
In particular, due to \rref{ta},
\beq \mbox{$T=\Tc=+\infty$ \quad if ~$\| u_0 \|_n \leqs \dd{\nu \over G_n}$}~; \label{coglob}
\feq
in this case $u$ is global and, if $\nu > 0$, it decays exponentially.
\end{prop}
\salto
\textbf{Other sufficient conditions for
global existence.}
Consider again the Cauchy problem \rref{cau} with
zero external forcing and any datum $u_0 \in \HM{n+2}$;
let $u \in C([0, T), \HM{n+2}) \cap C^1([0,T), \HM{n})$ be
the maximal solution.
An obvious implication of Proposition \ref{prozero} is the following.
\begin{prop}
\label{coro}
\textbf{Corollary.} Assume
\beq \| u(t_1) \|_n \leqs {\nu \over G_n}
\qquad \mbox{for some $t_1 \in [0,T)$}~. \label{assume} \feq
Then:
\beq T = + \infty, \qquad \| u(t) \|_n \leqs
{\| u(t_1) \|_n e^{-\nu (t-t_1)} \over 1 - G_n \| u(t_1) \|_n \, e_{\nu}(t-t_1)}
 \quad \mbox{for $t \in [t_1, +\infty)$}~,
\label{tes0} \feq
with $e_{\nu}$ as in Eq. \rref{enu}.
\end{prop}
\textbf{Proof.} The function $u \restriction [t_1, T)$
is the maximal solution of the Cauchy problem with datum $u(t_1)$ \textsl{specified
at time $t_1$}, rather than at time $0$. By Eqs. \rref{cobg} \rref{coglob},
with an obvious shift in time, we obtain the thesis \rref{tes0}. \fine
A consequence of the previous result, of more practical use, is
the following.
\begin{prop}
\label{coro2}
\textbf{Corollary.}
Let $\ua \in C([0, \Ta), \HM{n+2}) \cap C^1([0,\Ta), \HM{n})$ be any approximate
solution of \rref{cau},
with estimators $\ep_n, \delta_n, \Dd_n, \Dd_{n+1}$;
assume that the corresponding control inequalities
\rref{cont1} \rref{cont2} have a solution $\Rr_n \in C^1([0,\Tc),[0,+\infty))$,
with $\Tc \in (0,\Ta]$. Finally, assume that
\beq (\Dd_n + \Rr_n)(t_1) \leqs {\nu \over G_n}
\qquad \mbox{for some $t_1 \in [0,\Tc)$}~. \label{iprn} \feq
Then:
\beq T = + \infty, \quad \| u(t) \|_n \leqs
{(\Dd_n + \Rr_n)(t_1) e^{-\nu (t-t_1)} \over 1 - G_n (\Dd_n + \Rr_n)(t_1) \, e_{\nu}(t-t_1)}
~\mbox{for $t \in [t_1, +\infty)$}~.
\label{tes1} \feq
\end{prop}
\textbf{Proof.} By Proposition \ref{main}
we have $\| u(t) - \ua(t) \|_n \leqs \Rr_n(t)$
for all $t \in [0,\Tc)$ and, in particular, for $t=t_1$;
we further write
$\| u(t_1) \|_n \leqs \| \ua(t_1) \|_n +
\| u(t_1) - \ua(t_1) \|_n$, which implies
\beq \| u(t_1) \|_n \leqs (\Dd_n + \Rr_n)(t_1)~. \label{furth1} \feq
Eq. \rref{furth1} and the assumption \rref{iprn} gives the
inequality
$$ \| u(t_1) \|_n \leqs {\nu \over G_n}~, $$
which has the form \rref{assume}. By
the previous corollary, we have Eq. \rref{tes0} and
this result, combined with \rref{furth1},
gives the thesis \rref{tes1}. \fine
\salto
\textbf{A general result, under conditions of exponential decay
($\boma{\nu>0}$) or boundedness ($\boma{\nu=0}$) for the approximate solution.}
In this paragraph we exhibit a solution of the control
equations \rref{econt1} \rref{econt2}, holding
for \textsl{any} approximate solution whose estimators have certain
features. Such features are described
in the forthcoming Eq. \rref{dec1};
these indicate that the norms of orders $n$, $n+1$ of $\ua$
behave like $e^{-\nu t}$, while the $n$-th norm of the differential
error behaves like $e^{-2 \nu t}$. (In the NS case
these are conditions of exponential decay, while in the Euler
case they simply indicate the boundedness
of $\ua$ and its error.) \par
The assumption that $\ua$ behaves like $e^{-\nu t}$ corresponds
to a rather typical behavior of the approximate solutions
under zero external forcing; for example, this behavior occurs
in the case of the Galerkin approximate solutions
discussed in the next section (see Lemma \ref{rem64} and the subsequent
Remark). The differential error
of $\ua$ typically behaves like $e^{-2 \nu t}$ when
$\ua$ is bounded by $e^{-\nu t}$ and the differential error
depends only on the quadratic function $\PPP(\ua,\ua)$; again,
this situation occurs in the example of the Galerkin approximate
solutions with no forcing (see Lemma \ref{lem67}). \par
Other approximation methods suggested for the NS equations
($\nu > 0$) yield approximate solutions with a behavior as in
\rref{dec1}; as an example, this happens
if one assumes no external forcing
and takes for $\ua$ the truncation to any order of the power series
solution introduced by Sinai \cite{Sin}.
\begin{prop}
\label{lemdec}
\textbf{Lemma.} For any $\nu \in [0,+\infty)$, consider for \rref{cau}
an approximate solution
$\ua \in C([0, +\infty), \HM{n+2}) \cap C^1([0,+\infty), \HM{n})$.
Assume this to have growth and differential error estimators of the forms
\beq \Dd_n(t) := D_n e^{-\nu t}~, \quad \Dd_{n+1}(t) := D_{n+1} e^{-\nu t}~,
\quad \EE_n(t) := E_n e^{-2 \nu t}~,
\label{dec1}
\feq
with $D_n, D_{n+1} \in (0,+\infty)$, $E_n \in [0,+\infty)$; furthermore, assume this
to have any datum error estimator $\delta_n \in [0,+\infty)$.
From the above constants, let us define
\beq \De_n := {1 \over 2} (G_n D_n + K_{n} D_{n+1})~; \label{deden} \feq
furthermore, assume the ``error bound''
\beq G_n E_n < \De^2_n \label{erbound} \feq
and define
\beq W^{\pm}_n := \De_n \pm \sqrt{\De^2_n - G_n E_n}~. \label{dewn} \feq
Finally, consider the function $t \in [0,+\infty) \mapsto e_{\nu}(t)$
defined by \rref{enu}, and put
\beq \et_n : [0,+\infty) \vain [1,+\infty)~, \qquad
\et_n(t) := e^{\dd{(W^{+}_n - W^{-}_n) e_{\nu}(t)}}~. \label{etn} \feq
 Then, (i) (ii) hold. \par \noindent
(i) The control equations \rref{econt1} \rref{econt2} with the above estimators take the form
\beq {d \Rr_n \over d t} = - \nu \Rr_n + 2 F_n \, e^{-\nu t} \Rr_n +
G_{n} \Rr_n^2 + E_n e^{-2 \nu t}~, \qquad \Rr_n(0) = \delta_n~.
\label{esticc} \feq
(ii) The Cauchy problem \rref{esticc} has the solution $\Rr_n \in C^1([0,\Tc), [0,+\infty))$,
determined as follows:
\beq \Tc := \left\{ \barray{ll}
+ \infty & \mbox{if $\nu > 0$, ${W^{+}_n + G_n \delta_n \over W^{-}_n + G_n \delta_n} \geqs
e^{(W^{+} - W^{-})/\nu}$}~, \\
- \dd{1 \over \nu} \log \left(1 - \dd{\nu \over W^{+}_n -W^{-}_n}
\log {W^{+}_n + G_n \delta_n \over W^{-}_n + G_n \delta_n} \right)
& \mbox{if $\nu > 0$, ${W^{+}_n + G_n \delta_n \over W^{-}_n + G_n \delta_n} <
e^{(W^{+} - W^{-})/\nu}$}~, \\
\dd{1\over W^{+}_n -W^{-}_n}
\log {W^{+}_n + G_n \delta_n \over W^{-}_n + G_n \delta_n}  & \mbox{if $\nu=0$~;}
\farray \right. \label{tb} \feq
\beq \Rr_n(t) := {1\over G_n} {W^{+}_n (W^{-}_n + G_n \delta_n) \et_n(t)
- W^{-}_n (W^{+}_n + G_n \delta_n) \over
(W^{+}_n + G_n \delta_n) - (W^{-}_n + G_n \delta_n) \et_{n}(t)} \,e^{-\nu t}\, \mbox{for
$t \in [0,\Tc)$}. \label{eqrn} \feq
\end{prop}
\textbf{Proof.} (i) Obvious. \par \noindent
(ii) We write the unknown solution
$\Rr_n$ of \rref{esticc} as
\beq \Rr_n(t) =  \Zz_n(t) \, e^{-\nu t} \feq
where $\Zz_n \in C^1([0,\Tc), [0,+\infty))$ (and $\Tc$) are to be found.
Eq. \rref{esticc} is equivalent to the Cauchy problem
\beq {d \Zz_n \over d t} =
(E_n + 2 \De_n \Zz_n + G_n \Zz_n^2) \, e^{-\nu t}~, \qquad \Zz_n(0) = \delta_n~,
\label{estizz} \feq
with $\De_n$ as in \rref{deden}; this is solved by the quadrature formula
\beq \int_{\delta_n}^{\Zz_n(t)} \!\!\! {d z \over
E_n + 2 \De_n z + G_n z^2} = \int_{0}^t \!\! d s  \, e^{-\nu s}~.
\label{estiqq} \feq
On the other hand
\beq \int_{0}^t d s \,e^{-\nu s} = e_{\nu}(t)~; \label{imp1} \feq
furthermore $E_n + 2 \De_n z + G_n z^2 = (1/G_n)(W^{+}_n + G_n z)(W^{-}_n + G_n z)$,
which implies
\beq \int_{\delta_n}^{\Zz_n(t)} {d z \over
E_n + 2 \De_n z + G_n z^2} = - {1 \over W^{+}_n - W^{-}_n}
\left[ \log {W^{+}_n + G_n z \over W^{-}_n + G_n z} \right]^{z=\Zz_n(t)}_{z=\delta_n}~.
\label{imp2} \feq
Inserting Eqs. \rref{imp1} \rref{imp2} into Eq. \rref{estiqq} we easily
obtain an explicit expression for $\Zz_n(t)$ which implies
Eq. \rref{eqrn} for $\Rr_n(t) = \Zz_n(t) e^{-\nu t}$. The interval
$[0,\Tc)$ where $\Zz_n$ is well defined is also
easily determined, by the same manipulation employed to make explicit
$\Zz_n(t)$. \fine
From the previous lemma and from the main theorem
on approximate solutions (Proposition \ref{main}), we obtain the following result.
\begin{prop}
\label{prodec}
\textbf{Proposition.}
Consider the Euler/NS Cauchy problem \rref{cau},
and an approximate solution $\ua \in C([0, +\infty), \HM{n+2}) \cap C^1([0,+\infty), \HM{n})$;
assume this to have growth and differential error estimators
fulfilling all conditions in Lemma \ref{lemdec}, and
define $\Tc$, $\Rr_n$ via Eqs. \rref{tb} \rref{eqrn} of the same lemma.
Now, consider the maximal solution
$u \in C([0, T), \HM{n+2}) \cap C^1([0,T), \HM{n})$ of \rref{cau};
then
\beq T \geqs \Tc~, \qquad \| u(t) - \ua(t) \| \leqs \Rr_n(t) \qquad \mbox{for $t \in [0,\Tc)$}~.
\feq
(In particular, $u$ is global under the conditions
giving $\Tc = + \infty$ in Eq. \rref{tb}.)
\end{prop}
\begin{rema} \textbf{Remark.}
One could easily derive a variant of Lemma \ref{lemdec}
(and Proposition \ref{prodec}) dealing with the limit case $G_n E_n = \De^2_n$,
where $W^{-}_n = W^{+}_n$. The expressions for $\Tc$ and $\Rr_n$
in this case coincide with the ones derivable
from Eqs. \rref{tb} \rref{eqrn} in the limit $W^{-}_n \vain W^{+}_n$.
The results in Lemma \ref{lemzero} (and Proposition \ref{prozero}) could also be derived
from Lemma \ref{lemdec} (and Proposition \ref{prodec}), in the limit case $E_n=0$ and $D_{n}, D_{n+1} \vain 0$.
\end{rema}
\section{The Galerkin approximate solutions for
the Euler/NS equations, and their errors}
\label{galegen}
Throughout this section, we consider a set $G$ with the following features:
\beq \G \subset \Zd_0~, \qquad \G~\mbox{finite}~, \qquad
k \in \G \Leftrightarrow - k \in \G~. \label{modesg} \feq
Hereafter we write $\prec e_k \succ_{k \in G}$ for the linear subspace of $\Dz$
made of the sums $\sum_{k \in G} v_k e_k$ (which is, in fact, contained
in $\bb{C}^{\infty}$).
The forthcoming definitions follow closely
the setting proposed in our previous work \cite{due} for the Galerkin
method.
\salto
\textbf{Galerkin subspaces and projections.} We define them as follows.
\begin{prop}
\textbf{Definition.} The Galerkin
subspace and projection corresponding to $\G$ are, respectively:
\beq \HG := \Dsz \cap \prec e_k \succ_{k \in G} =
\{ \sum_{k \in \G} v_k e_k~|~v_k \in \complessi^d, \overline{v_k} = v_{-k}, k
\sc \, v_k = 0~\mbox{for all $k$} \}~; \label{hg} \feq
\beq \EG : \Dsz~ \vain \HG~, \qquad v =
\sum_{k \in \Zd_0} v_k e_k \mapsto \EG v := \sum_{k \in \G} v_k e_k~. \label{pg} \feq
\ffine
\end{prop}
It is clear that
\beq \HG \subset \bb{C}^{\infty} \cap \bb{\DD}'_{\so};~
\Delta(\HG) = \HG;~
\HG \subset \HM{m},~\EG(\HM{m}) = \HG~\mbox{for all $m \in \reali$}~. \feq
The following result is useful in the sequel.
\begin{prop}
\textbf{Lemma.} Let $m, p \in \reali$, $m \leqs p$ and $v \in
\HM{p}$. Then, \beq \| (1 - \EG) v \|_m \leqs {\| v \|_{p} \over
|\G|^{p-m}}~, \qquad |\G| := \min_{k \in \Zd_0 \setminus \G}
|k|~.\label{egmp} \feq
\end{prop}
\textbf{Proof.} We have $(1 - \EG) v = \sum_{k \in \Zd_0 \setminus G} v_k e_k$; thus,
\beq \| (1 - \EG) v \|^2_m = \sum_{k \in \Zd_0 \setminus G} |k|^{2 m} | v_k |^2 =
\sum_{k \in \Zd_0 \setminus G} {|k|^{2 p} \over | k |^{2 p- 2 m}} | v_k |^2  \feq
$$ \leqs \Big(\sup_{k \in \Zd_0 \setminus G} {1 \over |k|^{2 p-2 m}} \Big)
\sum_{k \in \Zd_0 \setminus G} |k|^{2 p} | v_k |^2  \leqs
{1 \over |\G|^{2 p-2 m}} \, \| v \|^2_{p}~. $$
\fine
\salto
\textbf{Galerkin approximate solutions.}
Let us be given
$\nu \in [0,+\infty)$, $f \in C([0,+\infty), \Dsz)$ and $u_0 \in \Dsz$.
\begin{prop}
\textbf{Definition.} The \textsl{Galerkin
approximate solution} of the Euler/NS equations corresponding to
the datum $u_0$, to the forcing $f$ and to the set of modes $G$ is the maximal
(i.e., nonextendable) solution
$u_{f, u_0, G} \equiv \ug$
of the following Cauchy problem, in the finite dimensional space $\HG$:
\beq \mbox{Find $\ug \in C^1([0,\Tg), \HG)$ such that} \label{galer} \feq
$$ {d \ug \over dt} = \nu \Delta \ug + \EG \PPP(\ug, \ug) +
\EG f~, \qquad \ug(0) = \EG u_0~. $$
\ffine
\end{prop}
Eq. \rref{galer} describes the Cauchy problem for an ODE in the
finite-dimensional vector space $\HG$; this relies
on the continuous function
$\HG \times [0,+\infty) \vain \HG$, $(v, t) \mapsto \nu \Delta v + \EG \PPP(v, v) +
\EG f(t)$, which is smooth with respect to the variable $v$. So, the standard
theory of ODEs in finite dimensional spaces grants existence
and uniqueness for the solution of \rref{galer}. \par
The following facts are known (see, e.g., \cite{Tem});
a proof suitable for the present setting is given, for completeness,
in Appendix \ref{appegal}.
\begin{prop}
\label{rem64}
\textbf{Lemma.} (i) For any
$\nu \in [0,+\infty)$, $f \in C([0,+\infty), \Dsz)$ and $u_0 \in \Dsz$,
the maximal solution of problem \rref{galer} is in fact global:
\beq \Tg = + \infty~. \label{tglob} \feq
(ii) If the forcing $f$ is identically zero, one has
\beq \| \ug(t) \|_{L^2}
\left\{ \barray{ll} = \| \EG u_0 \|_{L^2} & \mbox{if $\nu=0$, $t \in [0,+\infty)$,} \\
\leqs \| \EG u_0 \|_{L^2} \, e^{-\nu t} & \mbox{if $\nu>0$, $t \in [0,+\infty)$.} \farray \right.
\label{stacit} \feq
(iii) Let $\nu \in [0,+\infty)$. For any $f \in C([0,+\infty), \Dsz)$
one has
\beq \| \ug(t) \|_{L^2} \leqs \Big(\| \EG u_0 \|_{L^2} + \int_{0}^t d s \, e^{\nu s}
\| \EG f(s) \|_{L^2} \Big) \, e^{-\nu t} \quad \mbox{for $t \in [0,+\infty)$}~. \label{comb1} \feq
In particular, let
\beq
J := \int_{0}^{+\infty} \!\!\!\! ds \, e^{\nu s} \| \EG f(s) \|_{L^2}  < + \infty~;
\label{comb2} \feq
then,
\beq \| \ug(t) \|_{L^2} \leqs \Big(\| \EG u_0 \|_{L^2} + J \Big) \, e^{-\nu t}
\quad \mbox{for $t \in [0,+\infty)$}~. \label{comb3} \feq
\end{prop}
\begin{rema}
\textbf{Remark.}
Let $\nu \in [0,+\infty)$. If $f=0$ or, more generally, if $J < +\infty$,
the previous lemma gives a bound of the type $\| \ug(t) \|_{L^2} \leqs$ const. $e^{-\nu t}$.
This fact, with the equivalence of all norms on the finite
dimensional space $\HG$, implies the following: for any real $m$,
\beq \| \ug(t) \|_m \leqs U_m e^{-\nu t} \quad \mbox{for $t \in [0,+\infty)$}~, \label{teorest}
\feq
with a suitable constant $U_m \in [0,+\infty)$ (depending on the initial datum
and on the forcing).
\end{rema}
\begin{prop}
\textbf{Definition.}
(i) For all $k \in \Zd_0$,
$f_k \in C([0,+\infty), \complessi^d)$ and $u_{0 k} \in \complessi^d$
are the Fourier components
of the forcing and the initial datum: $f(t) = \sum_{k \in \Zd_0} f_k(t) e_k$,
$u_0 = \sum_{k \in \Zd_0} u_{0 k} e_k$. \par \noindent
(ii) For $k \in G$, $\ga_{G k} \equiv \ga_k
\in C^1([0,+ \infty), \complessi^d)$ are the Fourier components of $\ug$:
\beq \ug(t) = \sum_{k \in G} \ga_{k}(t) e_k ~.\feq
(Note the relations $\overline{f_k} = f_{-k}$,
$k \sc f_{k}=0$, and the analogous relations for $u_{0 k}$, $\ga_k$.)
\ffine
\end{prop}
Let us review a well known fact.
\begin{prop}
\textbf{Proposition.} Under the correspondence $u_G \mapsto (\gamma_k)$,
Eq. \rref{galer} for $\ug$ is equivalent
to the following problem: find a family of functions
$\ga_k \in C^1([0,+ \infty), \complessi^d)$ ($k \in G$) such that
\beq {d \ga_k \over d t} = - \nu |k|^2 \ga_k - {i \over (2 \pi)^{d/2}}
\sum_{h \in G} [ \ga_h \sc (k-h) ] \LP_k \ga_{k-h} + f_k~,
\qquad \ga_k(0) = u_{0 k} \label{galecomp} \feq
(intending $\ga_{k-h}(t) := 0$ if $k - h \not \in G$; recall that $\LP_k$
is the orthogonal projection of $\complessi^d$ onto $k^{\perp}$).
A family of functions $\gamma_k$ ($k \in G$) fulfilling Eqs.
\rref{galecomp} automatically possesses the properties
$\overline{\gamma_{k}} = \ga_{-k}$ and $k \sc \ga_k=0$.
\end{prop}
\textbf{Proof.} Clearly, Eq. \rref{galer} is equivalent to
$$ \left({d \ug \over dt}\right)_k = \nu (\Delta \ug)_k + \PPP(\ug, \ug)_k +
f_k~, \qquad (\ug)_k(0) = u_{0 k} \qquad \mbox{for $k \in G$}~; $$
making explicit the above Fourier components via Eq. \rref{repp}, we obtain Eqs. \rref{galecomp}.
\parn
Let us show that Eqs. \rref{galecomp} imply
$\overline{\gamma_{k}} = \ga_{-k}$ and $k \sc \ga_k=0$. The first implication
is proved noting that the functions $(\overline{\gamma_{k}})_{k \in G}$
and $(\ga_{-k})_{k \in G}$ are solutions of the same
Cauchy problem. The second implication follows noting that, for
each $k \in G$,
$(d/d t) (k \sc \ga_k) = - \nu |k|^2 (k \sc \ga_k)$ and $(k \sc \ga_k)(0)=0$;
the solution of this Cauchy problem is identically zero.
\fine
\textbf{The Galerkin solutions in the general framework
for approximate solutions.} From now on we stick
to the framework of the previous sections, i.e.: we consider
the Sobolev  spaces of orders $n,n+1, n+2$, for a fixed
$n \in (d/2+1, + \infty)$; as in \rref{vis}--\rref{in},
we choose a viscosity $\nu \in [0,+\infty)$, a forcing
$f \in C([0,+\infty), \HM{n})$ and an initial datum
$u_0 \in \HM{n+2}$; we consider the Euler/NS Cauchy problem
\rref{cau}.
Having fixed a set of modes $G$ as in \rref{modesg}, we
regard the corresponding Galerkin solution $\ug$ as an
approximate solution of the Euler/NS Cauchy problem
\rref{cau}; our aim is to apply the general theory of the previous
sections with $\ua = \ug$.
\par
The desired application requires to give growth
and error estimators for $\ug$.
Of course, we have the tautological growth estimators
\beq \Dd_m(t) := \| \ug(t) \|_m = \sqrt{ \sum_{k \in G} |k|^{2 m} |\ga_k(t)|^2 }~, \feq
to be used with $m=n$ or $m=n+1$; these are employed systematically
in the sequel. Let us pass to the errors of $\ug$ and their estimators.
\begin{prop}
\label{lem67}
\textbf{Lemma.}
(i) The Galerkin solution $\ug$ has the datum error
\beq \ug(0) - u_0 = -(1 - \EG) u_0 = - \sum_{k \in \Zd \setminus G} u_{0 k} e_k
\label{ovv1} \feq
and its tautological estimator
\beq \delta_n := \| \ug(0) - u_0 \|_{n} =
\sqrt{\sum_{k \in \Zd_0 \setminus G} |k|^{2 n} |u_{0 k}|^2}~. \label{ovv2} \feq
There is a rougher estimator
\beq \| \ug(0) - u_0 \|_{n} \leqs \delta'_{n p}~, \qquad \delta'_{n p} :=
{\| u_0 \|_{p} \over |G|^{p-n}}\,,
\label{ovv3} \feq
where $p$ is any real number such that $p \geqs n$, $u_0 \in \HM{p}$. \par \noindent
(ii) The differential error of $\ug$ is
\beq e(\ug) = - (1 - \EG) \PPP(\ug, \ug) - (1 - \EG) f~; \label{prov1} \feq
the Fourier representation of the above summands is
\beq (1 - \EG) f = \sum_{k \in \Zd_0 \setminus G} f_k e_k~; \label{prov2} \feq
\beq (1 - \EG) \PPP(\ug, \ug) = \sum_{k \in dG} p_k e_k~, \label{prov3} \feq
$$ dG := (G + G) \setminus (G \cup \{0\})~, \qquad p_k :=
- {i \over (2 \pi)^{d/2}}
\sum_{h \in G} [ \ga_h \sc (k-h) ] \LP_k \ga_{k-h}~. $$
(In the above: $G + G := \{p+q | ~p, q \in G\}$; $\setminus$ is the usual
set-theoretical difference; again, we intend $\ga_{k-h} := 0$
if $k - h \not\in G$.)\par \noindent
(iii) The above terms in $e(\ug)$ have norms
\beq \| (1 - \EG) f \|_n = \sqrt{\sum_{k \in \Zd_0 \setminus G} |k|^{2 n} |f_k|^2}~, \label{norfa} \feq
\beq \| (1 - \EG) \PPP(\ug, \ug) \|_n =
\sqrt{\sum_{k \in dG} |k|^{2 n} |p_k|^2}~; \label{nora} \feq
these admit the bounds
\beq \| (1 - \EG) f \|_n \leqs { \| f \|_{p} \over |G|^{p-n}}~, \label{norfb} \feq
\beq \| (1 - \EG) \PPP(\ug, \ug) \|_n \leqs
{K_{q} \over |G|^{q-n}} \| \ug \|_{q} \| \ug \|_{q+1}~\label{norb} \feq
where: $p$ is any real number such that $p \geqs n$ and $f \in C([0,+\infty), \HM{p})$;
$q$ is any real number such that $q \geqs n$; $K_q \in (0,+\infty)$
is such that $\| \PPP(v, w) \|_q \leqs K_q \| v \|_q \| w \|_{q+1}$
for all $v \in \HM{q}$, $w \in \HM{q+1}$ (of course
$\| \ug \|_{q} = \sqrt{\sum_{k \in G} |k|^{2 q} |\ga_k|^2}$,
and similarly for $\| \ug \|_{q+1}$).
Thus, a differential error estimator $\ep_n$ of order $n$ for $\ug$ is
obtained setting
\beq \ep_n := \mbox{( r.h.s. of \rref{nora} or \rref{norb} )} +
\mbox{( r.h.s. of \rref{norfa} or \rref{norfb}\,)}~. \label{estim} \feq
\end{prop}
\textbf{Proof.} (i) Eqs. \rref{ovv1} \rref{ovv2} are obvious. Eq.
\rref{ovv3} follows writing $\| \ug(0) - u_0 \|_{n}  =
\| (1 - \EG) u_0 \|_{n}$ and using the inequality \rref{egmp}. \par \noindent
(ii) We have
$$ e(\ug) = {d \ug \over d t} - \nu \Delta \ug - \PPP(\ug, \ug) - f =
\EG \PPP(\ug, \ug) + \EG f - \PPP(\ug, \ug) - f~, $$
where the first equality is just the definition of the differential
error, and the second follows from \rref{galer}; so, Eq. \rref{prov1}
is proved. \par
The Fourier representation \rref{prov2} for $(1 - \EG) f$ follows
immediately from the definition of $\EG$. To derive
the representation \rref{prov3} for $(1 - \EG) \PPP(\ug, \ug)$,
let us consider the Fourier components of $\ug$ that we indicate
in any case with $\ga_k$, intending $\ga_k := 0$
for $k \in \Zd \setminus G$. Eq. \rref{repp} gives
$$ \PPP(\ug, \ug) = \sum_{k \in \Zd_0} p_k e_k~, $$
where, for any $k$, $p_k$ is defined following Eq. \rref{prov3}.
On the other hand: for each $k \in \Zd_0$, $p_k$ is a sum over $h \in G$
containing terms of the form $\gamma_{k -h}$; if $k \not \in G+G$,
for all $h \in G$ we have $k - h \not \in G$ (since $k - h \in G$
would imply $k = (k - h) + h \in G+G$);
$k-h \not \in G$ implies $\gamma_{k-h} =0$. Due to the above considerations
we have $p_k = 0$ for $k \not \in G+G$, whence
\beq \PPP(\ug, \ug) = \sum_{k \in (G + G) \setminus \{0\}} p_k e_k~. \feq
Application of $1-\EG$ removes from the above sum the terms with $k \in G$;
thus
\beq (1 - \EG) \PPP(\ug, \ug) = \sum_{k \in (G + G) \setminus (G \cup \{0\})} p_k e_k~, \feq
in agreement with Eq. \rref{prov3}. \par \noindent
(iii) Eqs. \rref{norfa}  \rref{nora} are straightforward consequences
of the Fourier representations \rref{prov2} \rref{prov3}.
Eq. \rref{norfb} follows immediately from \rref{egmp}.
To infer Eq. \rref{norb}, we note that \rref{egmp} and
the inequality involving $K_q$ imply
$$ \| (1 - \EG) \PPP(\ug, \ug) \|_n \leqs {K_{q} \over |G|^{q-n}} \| \PPP(\ug, \ug) \|_q
\leqs {K_{q} \over |G|^{q-n}} \| \ug \|_{q} \| \ug \|_{q+1}~. $$
Finally, Eq. \rref{prov1} implies
\beq \| e(\ug) \| \leqs \| (1 - \EG) \PPP(\ug, \ug) \|_n + \| (1 - \EG) f \|_n~; \feq
binding the two summands in the right hand side via Eqs. \rref{nora} or \rref{norb},
\rref{norfa} or \rref{norfb} we obtain the estimator $\ep_n$ in \rref{estim}.
\fine
\begin{rema}
\textbf{Remarks.}
(i) The exact expression \rref{nora} for $\| (1 - \EG) \PPP(\ug, \ug)\|_n$
involves a sum over the \textsl{finite} set $dG$; for each $k \in dG$,
the term $p_k$ in the above sum is itself a finite sum, that can be
computed explicitly from the Fourier coefficients $\ga_k = \ga_k(t)$
of the Galerkin approximate solution (assuming these ones to be known,
say, from the numerical integration of the system \rref{galecomp}).
However, when $G$ is large the set $dG$ is typically very large, and
this can make  too expensive the computation of the sum
over $dG$. In these cases one can use for
$\| (1 - \EG) \PPP(\ug, \ug)\|_n$ the bound \rref{norb};
this requires computation of
$\| \ug \|_{q} = \sqrt{\sum_{k \in G} |k|^{2 q} |\ga_k|^2}$
and of $\| \ug \|_{q+1}$, both of them less expensive since these sums are
over $G$, rather than $dG$. \parn
(ii) We are aware that, if the Galerkin equations \rref{galecomp}
are solved numerically by some standard method for ODEs,
one does not obtain the exact solutions
$\ga_k(t)$ ($k \in G$) but, rather, some approximants whose
distance from the $\ga_k$'s should be estimated on the grounds
of the employed numerical scheme. In the application presented
in the next section, relying on a relatively small set $G$ of modes,
the intrinsic error in the numerical integration of \rref{galecomp} has been regarded
as negligible; the situation would be different using a much larger
set of modes.
\ffine
\end{rema}
\section{An application of the previous framework for the Galerkin method}
\label{caunecas}
\textbf{A preliminary.} In this section
we frequently report the results of computations
performed with MATHEMATICA. A formula like
$r= a. b c d e...~$ must be intended as follows: computation of
the real number $r$ via MATHEMATICA produces as an output $a.b c d e$,
followed by other digits not reported for brevity.
\salto
\textbf{Setting up the problem.}
Throughout this section, we work in space dimension
\beq d = 3 ~, \feq
with any viscosity
$\nu \in [0,+\infty)$. We consider the Euler/NS equations on $\Tt$ with no external forcing,
in the Sobolev framework of order $n=3$. So, the Cauchy problem
\rref{cau} takes the form
\parn
\vbox{
\beq \mbox{\textsl{Find}}~
u \in C([0, T), \HM{5}) \cap C^1([0,T), \HM{3}) \quad \mbox{\textsl{such that}} \label{caun} \feq
$$ {d u \over d t} = \nu \Delta u + \PPP(u,u)~, \qquad u(0) = u_0~. $$
}
The initial datum $u_0$ in $\HM{5}$ (in fact, in $\HM{m}$ for any real $m$)
is chosen of this form:
\par
\vbox{
\beq u_0 = \sum_{k = \pm a, \pm b, \pm c} u_{0 k} e_k~, \label{unec} \feq
$$ a := (1,1,0),~~ b := (1,0,1),~~ c := (0,1,1)~; $$
$$ u_{0, \pm a} := (2 \pi)^{3/2} (1,-1,0)~, \quad u_{0, \pm b}
:= (2 \pi)^{3/2} (1,0,-1)~, \quad u_{0, \pm c} := (2 \pi)^{3/2} (0,1,-1)~. $$
}
\noindent
For $\nu=0$, this initial datum has been considered by Behr, Ne$\check{\mbox{c}}$as
and Wu \cite{Nec} as the origin of a possible blow-up for the Euler equations.
More precisely, the authors of the cited work try a solution $u$ of the Euler equations
in the form of a power series $u(t) = \sum_{i=0}^{+\infty} u_i t^i$, where
the zero order term corresponds to the initial datum, and $u_1, u_2,... : \Tt \vain \reali^3$
are determined recursively. By a means of $C\scriptstyle{++}$ program,
the authors compute the terms $u_i$ for $i=1,...,35$, allowing
to construct the partial sums $u_N(t) := \sum_{i=0}^{N} u_i t^i$
up to $N=35$. A merely ``experimental'' analysis of these partial sums and
of their $\HM{3}$ norms brings the authors to conjecture
that the exact solution $u(t)$ of the Euler Cauchy problem blows up
in $\HM{3}$ for $t \vain \Tn$, for some $\Tn \in (0.32, 0.35)$.
\par \noindent
For subsequent use, we mention that Eq. \rref{unec}
implies
\beq \| u_0 \|_m = \sqrt{3 \pi^{3} 2^{m + 5}} \feq
for any real $m$. This formula gives, for example,
\beq \| u_0 \|_1 = 77.15...,~~ \| u_0 \|_2  = 109.1...,~~\label{u03} \feq
$$ \| u_0 \|_3 = 154.3....,~~ \| u_0 \|_4 = 218.2....,~~
\| u_0 \|_5 = 308.6...~. $$
\textbf{Introducing our approach.}
In this section we
propose a different approach to the Euler Cauchy problem of
\cite{Nec}, that we apply as well to the NS case; so, we consider
the problem \rref{caun} \rref{unec}, for arbitrary $\nu \in [0,+\infty)$.
\parn
We refer to the general setting developed in the present
paper, using the spaces
$\HM{n}$, $\HM{n +1}$, $\HM{n+2}$ with
\beq n = 3~. \feq
This requires, amongst else, the numerical values
of two constants $K_3$, $G_3$ fulfilling for $n=3$ the ``basic inequality''
\rref{basineq} and the ``Kato inequality'' \rref{katineq}.
Due to the computations in \cite{cok} \cite{cog},
these can be taken as follows:
\beq K_3 = 0.323~, \qquad G_3 = 0.438~. \label{k3g3} \feq
To illustrate our setting, we start with a very elementary
result.
\salto
\textbf{A simple sufficient condition for global existence (and
exponential decay).}
According to Proposition \ref{prozero}, the solution $u$ of the
NS Cauchy problem \rref{caun} \rref{unec} is global and exponentially decaying
for $t \vain + \infty$, if
\beq \nu \geqs G_3 \| u_0 \|_3 = 67.58... \label{glosimp} \feq
(in the last passage, we have used the numerical values in \rref{u03} \rref{k3g3}
for $G_3$ and $\| u_0 \|_3$)
({\footnote{For completeness, we mention
a criterion for NS global existence in $\HM{1}(\Tt)$,
discussed in \cite{accau}; this can be written as
$\nu \geqs \| u_0 \|_1/0.407$, where $u_0$
is an arbitrary initial datum in $\HM{1}$.
With the datum $u_0$ in \rref{unec} (and the value of $\| u_0 \|_1$
in \rref{u03}), we conclude that \cite{accau} ensures
global existence in $\HM{1}$ for $\nu \geqs 189.5...~$.
By a known regularity theorem about NS equations
\cite{Lem}, this result of global existence in $\HM{1}$ also implies global existence
in $\HM{3}$; however, the condition $\nu \geqs 189.5...$
arising from the $H^1$ setting of \cite{accau}
is manifestly weaker than the condition \rref{glosimp}.}}).
As shown hereafter, a more refined application
of our setting for approximate solutions significantly improves the above condition:
in fact, in the next pages, combining this setting with
the Galerkin method we will be able to infer global existence for $\nu \gtrsim 8$.
\salto
\textbf{Going on in our approach: Galerkin approximants.}
The idea developed in the sequel, both in the Euler case ($\nu=0$)
and in the NS case ($\nu >0$), is the following:
to compute numerically the Galerkin approximate
solution $\ug$ for a suitable set of modes $G$;
to construct for it error and growth estimators,
in the Sobolev norms of orders $3$ or $4$;
to solve numerically the control equations \rref{econt1} \rref{econt2} for
an unknown function $\Rr_3 : [0,\Tc) \vain [0,+\infty)$. After
finding the solution of this control problem, we can grant
on theoretical grounds that the solution $u$ of the Euler/NS
Cauchy problem \rref{caun} exists at least up to time
$\Tc$, and that $\| u(t) - \ug(t) \|_3 \leqs \Rr_3(t)$ for $t \in [0,\Tc)$.\par
Our computation has been performed using Mathematica on a PC, with the
relatively small set of $150$ modes
\par
\vbox{
\beq G := S \cup -S~; \qquad -S := \{ -k~|~k \in S \}~; \label{defg} \feq
$$ S :=
\{ (0, 0, 2), (0, 1, -3), (0, 1, 1), (0, 1, 3), (0, 2, 0), (0, 2, 2), (0, 3, -1),
(0, 3, 1), (0, 3, 3), $$
$$ (1, -3, -2), (1, -3, 0), (1, -3, 2), (1, -2, -3),
(1, -2, -1), (1, -2, 1), (1, -2, 3), $$
$$  (1, -1, -2), (1, -1, 2), (1, 0, -3),
(1, 0, 1), (1, 0, 3),
(1, 1, -2), (1, 1, 0),
(1, 1, 2), (1, 2, -3), $$
$$ (1, 2, -1), (1, 2, 1), (1, 2, 3), (1, 3, -2), (1, 3, 0),
(1, 3, 2), (2, -3, -3),
(2, -3, -1), (2, -3, 1), $$
$$ (2, -3, 3), (2, -2, -2), (2, -2, 2), (2, -1, -3),
(2, -1, -1), (2, -1, 1),
(2, -1, 3), (2, 0, 0), $$
$$ (2, 0, 2), (2, 1, -3), (2, 1, -1), (2, 1, 1), (2, 1, 3),
(2, 2, -2), (2, 2, 0),
(2, 3, -3), (2, 3, -1), $$
$$ (2, 3, 1), (2, 3, 3), (3, -3, -2), (3, -3, 2), (3, -2, -3),
(3, -2, -1),
(3, -2, 1), (3, -2, 3), $$
$$ (3, -1, -2), (3, -1, 0), (3, -1, 2), (3, 0, -1), (3, 0, 1),
(3, 0, 3), (3, 1, -2),
(3, 1, 0), (3, 1, 2), $$
$$ (3, 2, -3), (3, 2, -1), (3, 2, 1), (3, 2, 3), (3, 3, -2),
(3, 3, 0), (3, 3, 2)  \}~.
$$
}
\noindent
The results we present here are somehow provisional; we plan to attack the problem by more powerful
numerical tools in a future work, using for $G$ a larger set of modes.
\salto
\textbf{A sketch of the operations to be done.}
The list of such operations is the following: \par \noindent
(i) First of all, one chooses a value $\nu \in [0,+\infty)$ for the viscosity,
and a finite time interval $[0,\Tgg)$ for the numerical computation
of the Galerkin solution. \par \noindent
(ii) The Galerkin solution $\ug$ for the set of modes
$G$ in \rref{defg} and for the initial datum $u_0$ is found numerically
on the chosen time interval $[0, \Tgg)$. More precisely, one solves
numerically the system of differential equations \rref{galecomp}
for the Fourier components $(\gamma_k)_{k \in \G}$ of $\ug$,
with the initial conditions $u_{0 k}$ corresponding to Eq. \rref{unec}.
We recall (see Lemma \ref{rem64}) that the general theory of the Galerkin solutions
for zero external forcing ensures the global existence
of $\ug$ (and its exponential decay, if $\nu > 0$); thus,
from a theoretical viewpoint there is no obstruction to
the computation of $\ug$ on any finite interval $[0,\Tgg)$.
\par \noindent
(iii) We apply to $\ug$ on $[0,\Tgg)$
the framework of the present paper for the approximate
Euler/NS solutions using (we repeat it) the spaces
$\HM{n}$, $\HM{n+1}$ and $\HM{n+2}$ with $n = 3$.
\par \noindent
(iv) Some important characters in our approach are
the norms
\beq \Dd_n(t) := \| \ug(t) \|_{n} =
\sqrt{\sum_{k \in G} |k|^{2 n} |\ga_k(t)|^2}
\qquad (n=3,4)~, \label{d34} \feq
which can be obtained from the numerical values $\ga_k(t)$
of the Fourier components.
Further, we need a datum error
estimator $\delta_3$ and a differential error estimator $\ep_3$.
On the other hand the datum error is zero
in this case, since the initial condition
$u_0$ is in the subspace $\HG$ spanned by the chosen set \rref{defg} of modes; thus,
we can take
\beq \delta_3 = 0~.\feq
As for the differential error estimator,
we take the precise expression coming from
Eqs. \rref{nora} \rref{estim}
(taking into account that the forcing $f$ is
zero in this case); these yield the expression
\beq \ep_3(t) = \sqrt{\sum_{k \in dG} |k|^{6} |p_k(t)|^2} \label{nora3} \feq
where $dG := (G + G) \setminus (G \cup \{0\})$ and
$p_k(t) := - i (2 \pi)^{-d/2}
\sum_{h \in G} [ \ga_h(t) \sc (k-h) ] \LP_k \ga_{k-h}(t)$, as in Eq. \rref{prov3}.
The computation of $\ep_3(t)$ following Eq. \rref{nora3}
is rather expensive: the set $dG$ consists of $929$ modes
and, for each one of them, one must perform the nontrivial
computation of $p_k(t)$. (For these reasons,
a computation with a set of modes much larger than the $G$ in
\rref{defg} would suggest to replace this $\ep_3$
with a rougher estimator,
coming from Eqs. \rref{norb} \rref{estim}). \par \noindent
(v) Now we consider the control equations \rref{econt1} \rref{econt2},
taking the form
\beq {d \Rr_3 \over d t} = - \nu \Rr_3
+ (G_3 \Dd_3 + K_3 \Dd_{4}) \Rr_3 + G_3 \Rr^2_3 + \ep_3~,
\label{econt13} \feq
\beq \Rr_3(0) = 0~, \label{econt23} \feq
with $K_3, G_3$ as in Eq.\rref{k3g3}
and $\Dd_3(t), \Dd_4(t), \ep_3(t)$ as in Eqs. \rref{d34} \rref{nora3}.
The unknown is a function $\Rr_3 \in C^1([0, \Tc), [0,+\infty))$,
with $0 < \Tc \leqs \Tgg$; this is determined numerically (with a
package for self-adaptive integration, allowing to detect a possible
blow-up of $\Rr_3$; in this case, $\Tc$ is the blow-up time). \par
Once we have $\Rr_3$, the general theory allows to state the following. \par\noindent
(a) The maximal solution $u$ of the Euler/NS Cauchy problem
\rref{caun} \rref{unec} has a domain containing $[0,\Tc)$, and
\beq \| u(t) - \ug(t) \|_3 \leqs \Rr_3(t) \qquad \mbox{for $t \in [0,\Tc)$} \label{ubr3} \feq
(see Proposition \ref{main}). \parn
(b) Suppose $\nu > 0$, and
\beq (\Dd_3 + \Rr_3)(t_1) \leqs {\nu \over G_3}
\qquad \mbox{for some $t_1 \in [0,\Tc)$}~. \label{iprn3} \feq
Then, the solution $u$ of the NS Cauchy problem \rref{caun}
is global, and
\beq \| u(t) \|_3 \leqs
{(\Dd_3 + \Rr_3)(t_1) \, e^{-\nu (t-t_1)} \over 1 - G_3 (\Dd_3 + \Rr_3)(t_1) \, e_{\nu}(t-t_1)}
 \quad \mbox{for $t \in [t_1, +\infty)$}~,
\label{tes} \feq
with $e_{\nu}$ as in Eq. \rref{enu} (see Corollary \ref{coro2}).
\salto
Hereafter, we report the results obtained from the above scheme for some choices of $\nu$.
\salto
\textbf{Case $\boma{\nu=0}$ }. First of all, the Galerkin
equations \rref{galecomp} have been solved for the set of modes
\rref{defg} on a time interval of length $\Tgg = 2$.
(This required, approximately, $15$ seconds of CPU time on our machine.)
In Figures 1a, 1b and 1c we report, as examples, the graphs of
$|\gamma_k(t)|$ for $t \in [0,2)$, in the cases
$k=(1,1,0)$, $k=(0,0,2)$ and $k=(0,1,-3)$, respectively. (Of course,
we could consider many alternative choices, such as plotting
the norms $|\gamma_k(t)|$ for other modes, or
the real parts of the components $\gamma^r_{k}(t)$ ($r=1,2,3$),
or the imaginary parts of the same components.) As a supplementary
information, we mention that the graphs of $|\gamma_k(t)|$
are identical in the three cases $k = (1,1,0), (1,0,1), (0,1,1)$). \par
After computing numerically all the $\gamma_k(t)$ (for $k \in G$),
one obtains from Eqs. \rref{d34} and \rref{nora3}
the norms $\Dd_n(t) := \| \ug(t) \|_n$
($n=3,4$) and the differential error estimator $\ep_3(t)$, for $t \in [0,2)$;
Figures 1d and 1e report the graphs of $\Dd_3(t)$ and $\ep_3(t)$
(to give some more detail, we mention that the functions
$\Dd_3(t)$, $\Dd_4(t)$ and $\ep_3(t)$ have been computed for a set
of sample values of $t \in [0,2)$, and then interpolated
by means of the MATHEMATICA algorithms. The computation of
$\ep_3(t)$ at each sample value $t$ is rather expensive,
for it involves a sum on the large set $dG$; the CPU time
is about $15$ seconds for each $t$ and, for this reason,
we have used only $30$ sample values). \par
The final step is the numerical solution of the Cauchy problem
\rref{econt1} \rref{econt2} for the unknown function $\Rr_3(t)$
(a task performed almost instantaneously by our PC).
The MATHEMATICA self-adaptive routines for ODEs indicate
a divergence of $\Rr_3(t)$ for $t \vain \Tc$, with
$\Tc = 0.06666...$\,. Figure 1f gives the graph of
$\Rr_3(t)$ for $t \in [0, \Tc)$.
The conclusion of these computations is the following:
the solution $u$ of the Euler Cauchy problem \rref{caun}
\rref{unec} has a domain containing $[0, \Tc)$, and
its distance from $\ug$ is bounded by $\Rr_3$
on this interval.
\salto
\textbf{Some remarks on the $\boma{\nu=0}$ case.}
We have already mentioned that
the blow-up of $u(t)$ is conjectured in \cite{Nec} for $t \vain
\Tn \in (0.32, 0.35)$. The basis of this conjecture
is an experimental analysis of the first 35 terms
in the power series (in time) solving formally
the Cauchy problem; in principle, this analysis
does not even prove existence
of $u(t)$ at \textsl{any} specified time $t < \Tn$.
Our lower bound
$\Tc = 0.06666...$ for the interval of existence of
$u$ is about $1/5$ of the blow-up time $\Tn$
suggested by \cite{Nec}, but it relies
on an analytic theory for approximate solutions,
their errors, etc., summarized by Proposition
\ref{main}; in this sense, it is theoretically
grounded. \par
Perhaps, our method could give a sensibly larger lower bound
$\Tc$ on the time of existence, when
implemented with a set $G$ of Galerkin modes
much larger than \rref{defg}. Alternatively,
one could apply our theoretical framework
using as an approximate solution the
partial sum $u_N(t) := \sum_{i=1}^N
u_i t^i$, with $N=35$ as in \cite{Nec}
or with a larger $N$. Both tasks require much more expensive numerical
computations (to be done with more appropriate hardware and software):
we plan to do this elsewhere. One cannot
exclude that an attack with more powerful
devices could finally result into a theoretically
grounded lower bound $\Tc$ larger than the
suspected blow-up time $\Tn$ of \cite{Nec};
however, at present this is just a hope.
\salto
\textbf{Cases $\boma{\nu = 3}$ and $\boma{\nu=7}$}.
We use again the Galerkin solution $\ug$ with $G$ as in \rref{defg}. Due to
the positivity of $\nu$, all components $\gamma_{k}(t)$
of the Galerkin solution decay exponentially for $t \vain + \infty$
(recall Lemma \ref{rem64});
the same happens of the norms $\Dd_3(t)$, $\Dd_4(t)$ and of the differential error
estimator $\ep_3(t)$, which are
essential objects for our purposes. \par
The system \rref{galecomp} for the Galerkin components $\gamma_k(t)$ has been solved numerically on
a time interval of length $\Tgg=1$ (which
required a CPU time of about $15$ seconds for $\nu=3$,
and $25$ seconds for $\nu=7$). Subsequently, $\Dd_3(t)$, $\Dd_4(t)$
and $\ep_3(t)$ have been computed
from the components $\gamma_k(t)$ and Eqs. \rref{d34} \rref{nora3}
(indeed, some interpolation has been done as in the case $\nu=0$;
as in that case, for the computation of $\ep_3(t)$
we have used only 30 sample values of $t$ in $[0,1)$,
with a CPU time of about $15$ seconds for each one). \par
For both the above values of $\nu$, the final step has been
the (very fast) numerical solution of the Cauchy problem
\rref{econt1} \rref{econt2} for the unknown function $\Rr_3(t)$.
According to the MATHEMATICA routines for ODEs,
$\Rr_3(t)$ diverges for $t \vain \Tc$, where
$\Tc = 0.09025...$ for $\nu=3$, and $\Tc =  0.2386...$ for $\nu=7$.
We repeat that these results grant
existence on a domain $\supset [0, \Tc)$
for the solution $u$ of the NS Cauchy problem \rref{caun}
\rref{unec}, and the bound \rref{ubr3} on
this interval. \par
As examples, in Figures 2a-2f we
have reported some details on computations
for $\nu=7$. More precisely,
Figures 2a, 2b and 2c are the graphs of
$|\gamma_k(t)|$ for $t \in [0,1)$, in the cases
$k=(1,1,0)$, $k=(0,0,2)$ and $k=(0,1,-3)$, respectively.
(In fact, the graphs of $|\gamma_k(t)|$
for $k=(1,0,1)$ and $k=(0,1,1)$
are identical to the graph
of the case $k = (1,1,0)$.)
Figures 2d and 2e are the graphs of $\Dd_3(t)$ and $\ep_3(t)$,
for $t \in [0,1)$.
Figure 2f gives the graph of
$\Rr_3(t)$ that, as anticipated,
diverges for $t \vain \Tc = 0.2386...$\, .
\salto
\textbf{Case $\boma{\nu = 8}$}.
Again, we have  used
the Galerkin solution $\ug$ with $G$ as in \rref{defg}.
All the components $\gamma_{k}(t)$,
the norms $\Dd_3(t)$, $\Dd_4(t)$ and the differential error estimator $\ep_3(t)$
decay exponentially for $t \vain + \infty$.
The system \rref{galecomp}
for the Galerkin components $\gamma_k(t)$ has been solved numerically on
a time interval of length $\Tgg=1$ (which
required a CPU time of about $25$ seconds).
The forthcoming
Figures 3a,3b,3c are the graphs of
$|\gamma_k(t)|$ for $t \in [0,1)$,
$k=(1,1,0)$, $k=(0,0,2)$ and $k=(0,1,-3)$, respectively.
Subsequently, the norms $\Dd_n(t)$
($n=3,4$) and the error $\ep_3(t)$ have been computed
from the components $\gamma_k(t)$ and Eqs. \rref{d34} \rref{nora3}
(making some interpolations,
as in the previous cases).
Figures 2d and 2e are the graphs of $\Dd_3(t)$ and $\ep_3(t)$,
for $t \in [0,1)$.
\par
The final step has been
the (very fast) numerical solution of the Cauchy problem
\rref{econt1} \rref{econt2} for the unknown function $\Rr_3(t)$.
Differently from all the previous cases,
the numerical solution $\Rr_3(t)$ determined by MATHEMATICA
is defined on the whole interval $[0,1)$;
its graph is reported in Figure 3f which suggests,
via some extrapolation, that $\Rr_3(t)$
should be defined for all $t \in [0,+\infty)$,
with $\Rr_3(t) \vain 0^{+}$ for $t \vain + \infty$;
of course, this
would imply global existence for the solution $u$
of the NS Cauchy problem \rref{caun} \rref{unec}.
\par
However, global existence of $u$ can even be inferred
without extrapolating the behavior of $\Rr_3$
outside the interval $[0,1)$. In fact,
global existence is granted if the condition \rref{iprn3}
$(\Dd_3 + \Rr_3)(t_1) \leqs \nu/G_3$ holds
at any instant $t_1 > 0$; in the present case
$\nu/G_3 = 18.26... $, and the numerical computation
performed in the interval $[0,1)$ shows that \rref{iprn3}
is satisfied for any $t_1 \in [0.1567..., 1)$.
In conclusion, we can take for granted that we have
global existence for the solution $u$ of the NS Cauchy problem.
Of course, $\| u(t) - \ua(t) \|_3$
is bounded by the numerically computed function
$\Rr_3(t)$, for $t \in [0,1)$. After choosing
a $t_1 \in [0.1567, + \infty)$, we obtain as well a bound
of the form \rref{tes0} for $\| u(t) \|_3$, which also implies
$\| u(t) \|_3$ to vanish exponentially for $t \vain + \infty$.
\par
For example, let us choose $t_1 = 0.9$. We have $\Dd_3(0.9) =
8.580... \times 10^{-5}$, $\Rr_3(0.9) = 0.06100...$; with the already known value $G_3 = 0.438$,
Eq. \rref{tes0} gives
\beq \| u(t) \|_3 \leqs
{0.0614 \, e^{-8 (t-0.9)} \over 1 + 0.00335 \, e^{-8 (t-0.9)}} \leqs 0.0614 \, e^{-8 (t-0.9)}
\qquad \mbox{for $t \in [0.9, + \infty)$}~. \label{espbound} \feq
(Here the first inequality follows
directly from \rref{tes0}, recalling that $e_8(t) = (1 - e^{ - 8 t})/8$;
the second inequality is obvious.)
\salto
\textbf{Summary of the previous results, and final comments.}
Our method to treat the Galerkin approximant $\ug$, with
$G$ as in \rref{defg}, has given the following results
about the Euler/NS Cauchy problem \rref{caun} \rref{unec}. \par \noindent
a) $\nu=0,3,7$: we can grant existence of the solution $u$
of \rref{caun} \rref{unec} on an interval containing
$[0, \Tc)$, with $\Tc =0.06666, 0.09025, 0.2386$, respectively, for
these three choices of $\nu$. We have $\| u(t) - \ug(t) \|_3
\leqs \Rr_3(t)$ for $t \in [0, \Tc)$, where $\Rr_3$ is computed
numerically solving \rref{econt1} \rref{econt2} (for $\nu=0$ and
$\nu=7$, the graph of $\Rr_3$ is reported in
Figures 1f and 2f). \par \noindent
b) $\nu=8$: we can grant global existence for the
solution $u$ of \rref{caun} \rref{unec}. For $t \in [0,1)$
we have a bound $\| u(t) - \ug(t) \|_3 \leqs \Rr_3(t)$,
with $\Rr_3$ obtained again from the numerical solution of
\rref{econt1} \rref{econt2}; the graph of this
function is reported in Figure 3f. For $t \in [0.9, +
\infty)$ we have a bound of the form \rref{espbound}
on $\| u(t) \|_3$, decaying exponentially for large $t$.
\par
By extrapolation, we are led to conjecture that results
similar to (a) should be obtained for all $\nu \in [0,\nu_{cr})$,
while results similar to (b) should be obtained for
all $\nu \in [\nu_{cr}, + \infty)$, for some $\nu_{cr} \in (7,8)$.
(From a qualitative viewpoint, this is just the behavior
described by Lemma \ref{lemdec} on the analytical
solution of the control inequality.
\par \noindent
\vbox{
\noindent
However,
here we are applying the control equalities
with the tautological growth and error estimators, given
directly by the numerical solution of the Galerkin
equations; these are more precise than the
simple analytical estimators considered
in the cited lemma.)
}
\begin{figure}
\parbox{3in}{
\includegraphics[
height=2.0in,
width=2.8in
]%
{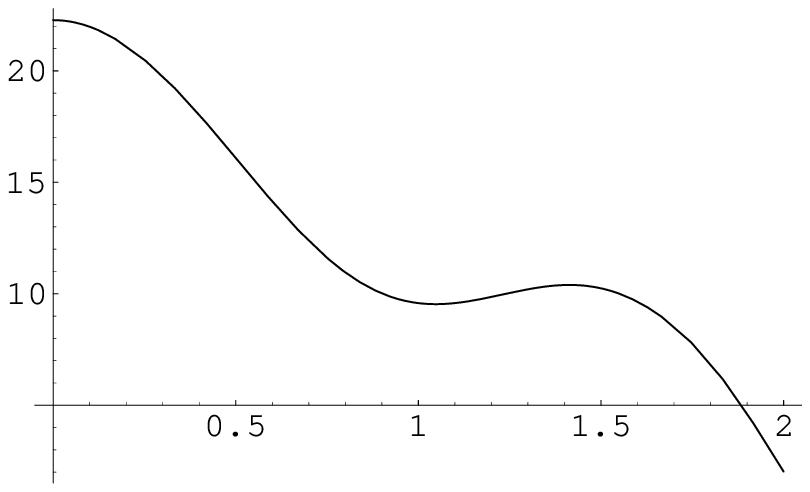}%
\par
{\footnotesize{
{\textbf{Figure 1a.~} $\nu=0$. Graph
of $|\gamma_{k}(t)|$ when $k=(1,1,0)$,
for $t \in [0,2)$.
}}
\par}
\label{f1a}
}
\hskip 0.4cm
\parbox{3in}{
\includegraphics[
height=2.0in,
width=2.8in
]%
{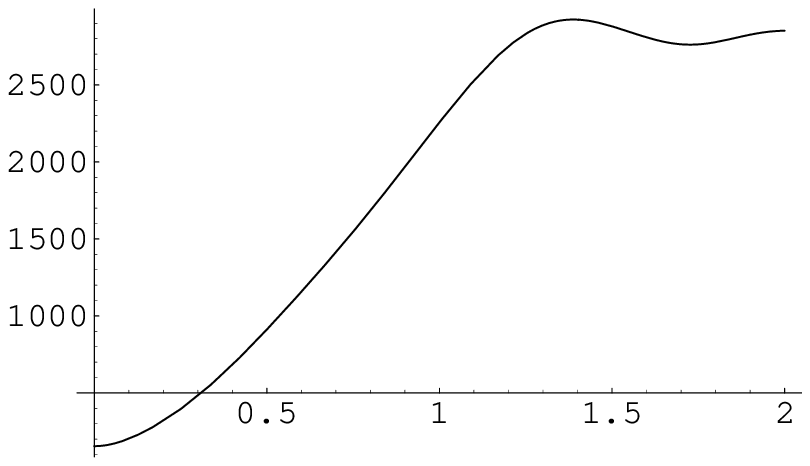}%
\par
{\footnotesize{
{\textbf{Figure 1d.~} $\nu=0$.
Graph of $\Dd_3(t) = \| \ug(t) \|_3$,
for $t \in [0,2)$.
One has $\Dd_3(0) = 154.3...$, $\Dd_3(0.02) = 156.4...$,
$\Dd_3(0.04) = 162.7...$, $\Dd_3(0.06) = 172.7...$~.
}}
\par}
\label{f1d}
}
\parbox{3in}{
\includegraphics[
height=2.0in,
width=2.8in
]%
{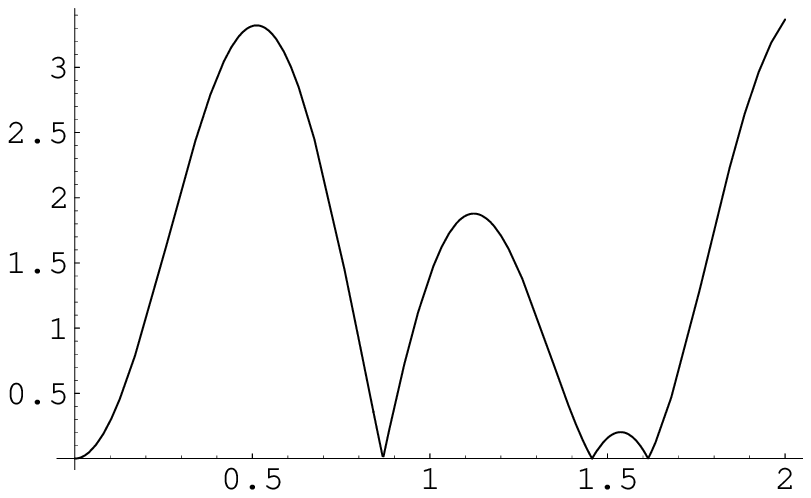}%
\par
{\footnotesize{
{\textbf{Figure 1b.~}
$\nu=0$. Graph
of $|\gamma_{k}(t)|$ when $k=(0,0,2)$,
for $t \in [0,2)$.
}}
\par}
\label{f1b}
}
\hskip 0.4cm
\parbox{3in}{
\includegraphics[
height=2.0in,
width=2.8in
]%
{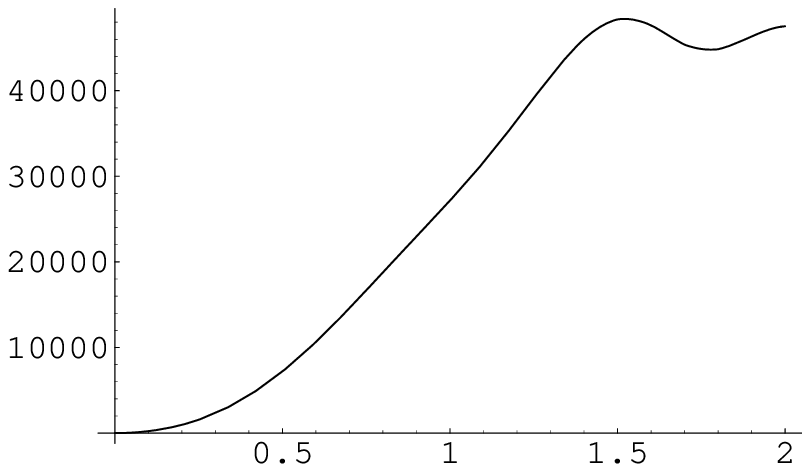}%
\par
{\footnotesize{
{\textbf{Figure 1e.~} $\nu=0$. Graph
of $\ep_3(t)$, for $t \in [0,2)$.
One has $\ep_3(0.02) = 9.027...$,
$\ep_3(0.04) = 36.17...$,
$\ep_3(0.06) = 81.65...$\, .
}}
\par}
\label{f1e}
}
\parbox{3in}{
\includegraphics[
height=2.0in,
width=2.8in
]%
{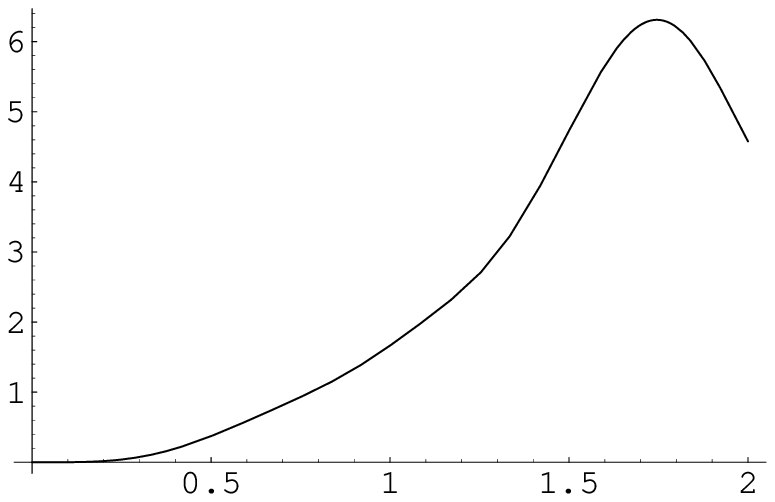}%
\par
{\footnotesize{
{\textbf{Figure 1c.~}
$\nu=0$. Graph
of $|\gamma_{k}(t)|$ when $k=(0,1,-3)$,
for $t \in [0,2)$.
}}
\par}
\label{f1c}
}
\hskip 0.4cm
\parbox{3in}{
\includegraphics[
height=2.0in,
width=2.8in
]%
{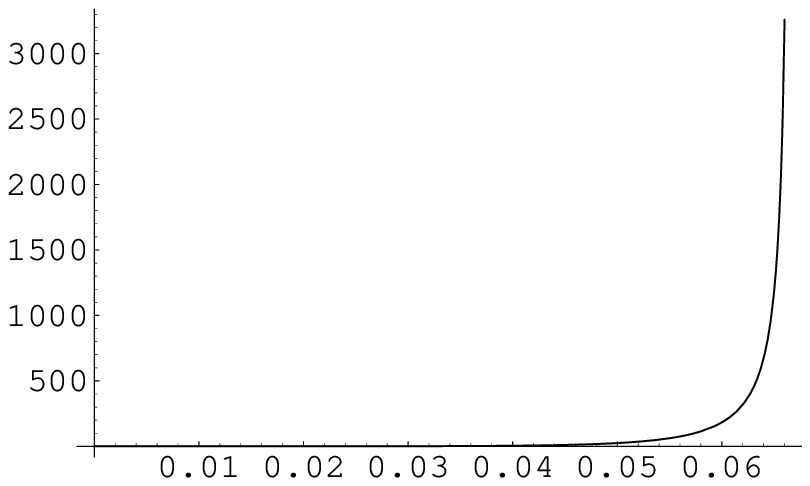}%
\par
{\footnotesize{
{\textbf{Figure 1f.~} $\nu=0$. Graph
of $\Rr_3(t)$; this function diverges
for $t \vain \Tc = 0.06666...$\,.
One has $\Rr_3(0) = 0$, $\Rr_3(0.02) = 0.1439...$,
$\Rr_3(0.04) = 4.685...$, $\Rr_3(0.06) = 182.3...$~.
}}
\par}
\label{f1f}
}
\end{figure}
\begin{figure}
\parbox{3in}{
\includegraphics[
height=2.0in,
width=2.8in
]%
{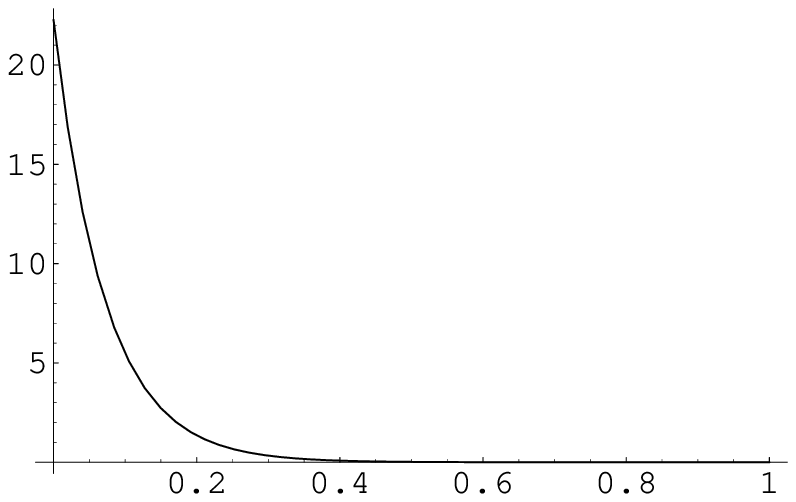}%
\par
{\footnotesize{
{\textbf{Figure 2a.~} $\nu=7$. Graph
of $|\gamma_{k}(t)|$ when $k=(1,1,0)$,
for $t \in [0,1)$.
}}
\par}
\label{f2a}
}
\hskip 0.4cm
\parbox{3in}{
\includegraphics[
height=2.0in,
width=2.8in
]%
{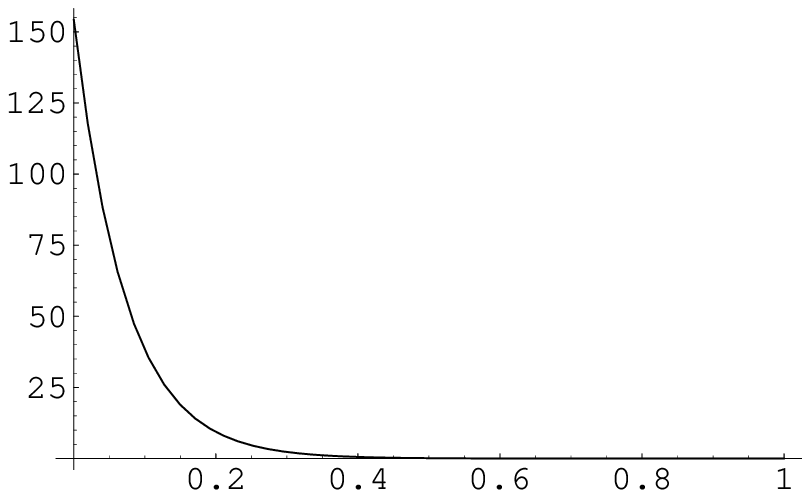}%
\par
{\footnotesize{
{\textbf{Figure 2d.~} $\nu=7$.
Graph of $\Dd_3(t)$,
for $t \in [0,1)$.
One has $\Dd_3(0) = 154.3...$, $\Dd_3(0.025) = 109.5...$, $\Dd_3(0.07) = 58.40...$,
$\Dd_3(0.15) = 18.89...$, $\Dd_3(0.23) = 6.153...$~.
}}
\par}
\label{f2d}
}
\parbox{3in}{
\includegraphics[
height=2.0in,
width=2.8in
]%
{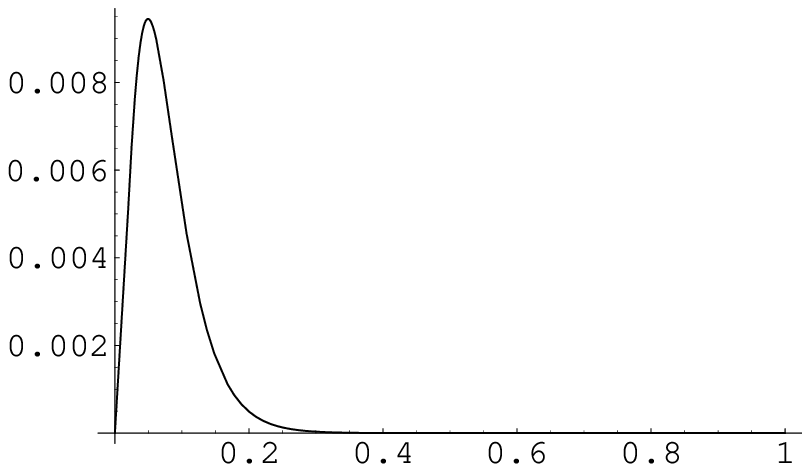}%
\par
{\footnotesize{
{\textbf{Figure 2b.~}
$\nu=7$. Graph
of $|\gamma_{k}(t)|$ when $k=(0,0,2)$,
for $t \in [0,1)$.
}}
\par}
\label{f2b}
}
\hskip 0.4cm
\parbox{3in}{
\includegraphics[
height=2.0in,
width=2.8in
]%
{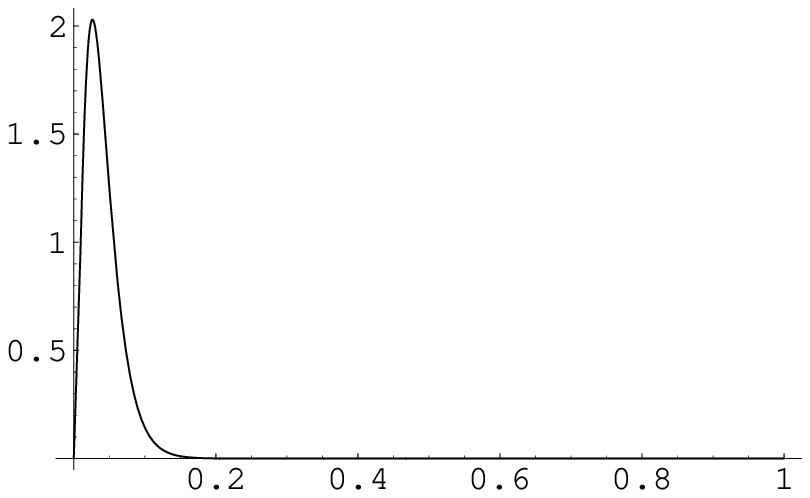}%
\par
{\footnotesize{
{\textbf{Figure 2e.~} $\nu=7$. Graph
of $\ep_3(t)$, for $t \in [0,1)$.
One has $\ep_3(0)=0$,
$\ep_3(0.025)=  2.024...$,
$\ep_3(0.07) = 0.5788...$,
$\ep_3(0.15) = 0.01089...$,
$\ep_3(0.23) = 1.338... \times 10^{-4}$.
}}
\par}
\label{f2e}
}
\parbox{3in}{
\includegraphics[
height=2.0in,
width=2.8in
]%
{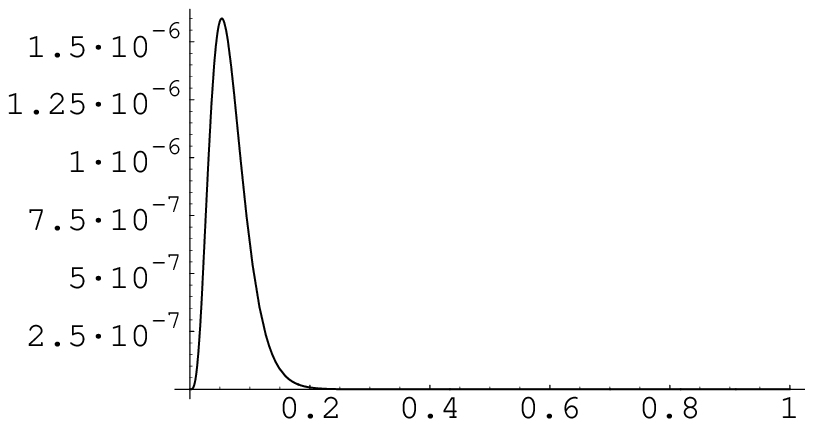}%
\par
{\footnotesize{
{\textbf{Figure 2c.~}
$\nu=7$. Graph
of $|\gamma_{k}(t)|$ when $k=(0,1,-3)$,
for $t \in [0,1)$.
}}
\par}
\label{f2c}
}
\hskip 0.4cm
\parbox{3in}{
\includegraphics[
height=2.0in,
width=2.8in
]%
{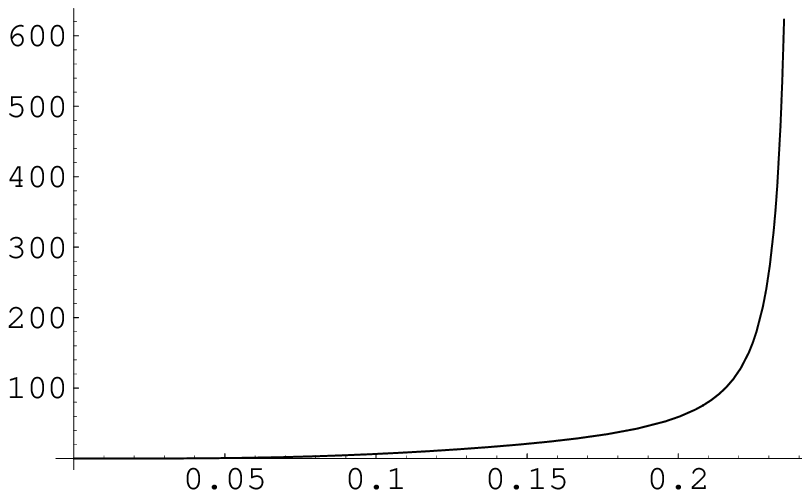}%
\par
{\footnotesize{
{\textbf{Figure 2f.~} $\nu=7$. Graph
of $\Rr_3(t)$; this function diverges
for $t \vain \Tc = 0.2386...$\,.
One has $\Rr_3(0) = 0$, $\Rr_3(0.07) = 2.096...$,
$\Rr_3(0.15) = 20.90...$, $\Rr_3(0.23) = 265.0...$~.
}}
\par}
\label{f2f}
}
\end{figure}
\begin{figure}
\parbox{3in}{
\includegraphics[
height=2.0in,
width=2.8in
]%
{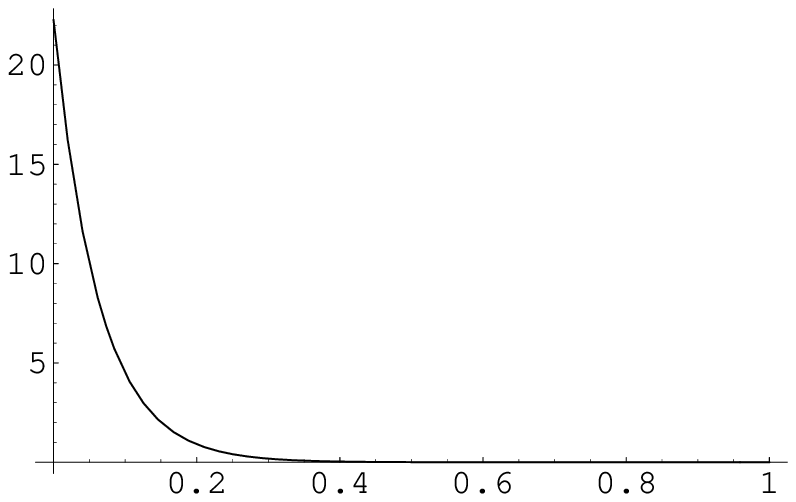}%
\par
{\footnotesize{
{\textbf{Figure 3a.~} $\nu=8$. Graph
of $|\gamma_{k}(t)|$ when $k =(1,1,0)$,
for $t \in [0,1)$.
}}
\par}
\label{f3a}
}
\hskip 0.4cm
\parbox{3in}{
\includegraphics[
height=2.0in,
width=2.8in
]%
{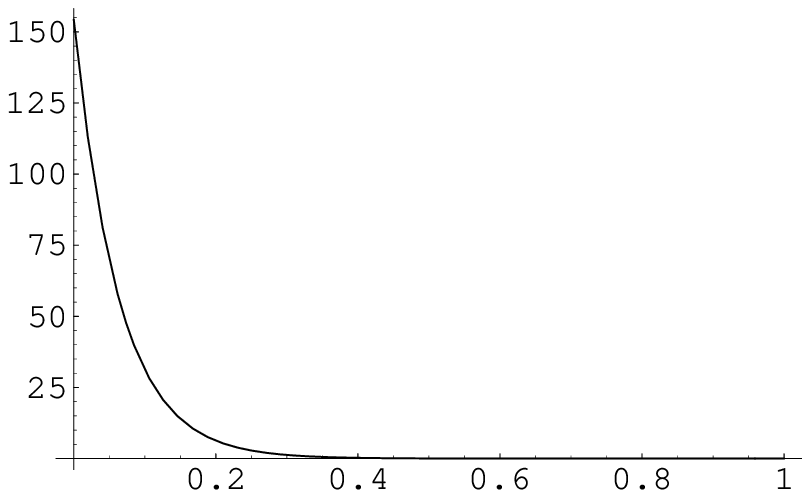}%
\par
{\footnotesize{
{\textbf{Figure 3d.~} $\nu=8$.
Graph of $\Dd_3(t)$,
for $t \in [0,1)$.
One has $\Dd_3(0) = 154.3...$, $\Dd_3(0.2) = 6.280...$,
$\Dd_3(0.4) = 0.2559...$,
$\Dd_3(0.6) = 0.01043...$, $\Dd_3(0.9) = 8.580... \times 10^{-5}...$~.
}}
\par}
\label{f3d}
}
\parbox{3in}{
\includegraphics[
height=2.0in,
width=2.8in
]%
{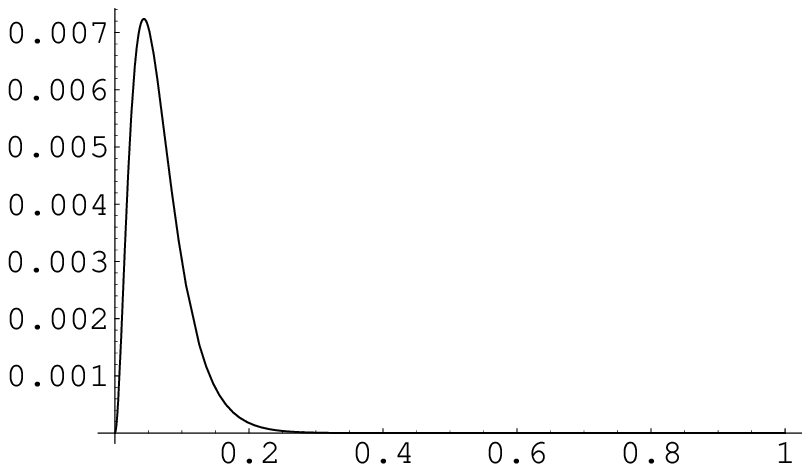}%
\par
{\footnotesize{
{\textbf{Figure 3b.~}
$\nu=8$. Graph
of $|\gamma_{k}(t)|$ when $k=(0,0,2)$,
for $t \in [0,1)$.
}}
\par}
\label{f3b}
}
\hskip 0.4cm
\parbox{3in}{
\includegraphics[
height=2.0in,
width=2.8in
]%
{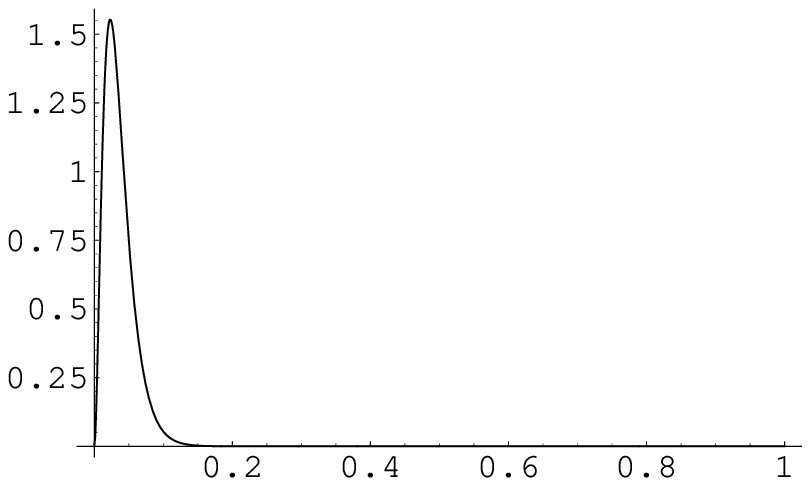}%
\par
{\footnotesize{
{\textbf{Figure 3e.~} $\nu=8$. Graph
of $\ep_3(t)$, for $t \in [0,1)$.
One has $\ep_3(0)=0$,
$\ep_3(0.023)= 1.553...$,
$\ep_3(0.07) = 0.2844...$,
$\ep_3(0.15) = 0.002638...$,
$\ep_3(0.23) = 1.593... \times 10^{-5}$.
}}
\par}
\label{f3e}
}
\parbox{3in}{
\includegraphics[
height=2.0in,
width=2.8in
]%
{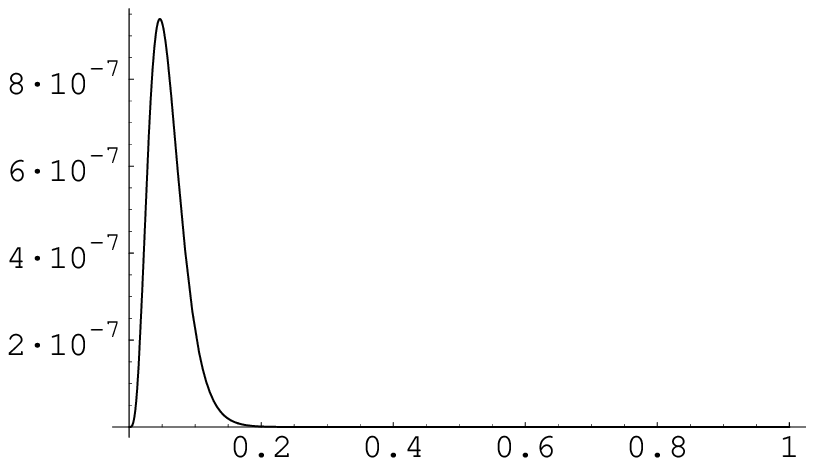}%
\par
{\footnotesize{
{\textbf{Figure 3c.~}
$\nu=8$. Graph
of $|\gamma_{k}(t)|$ when $k=(0,1,-3)$,
for $t \in [0,1)$.
}}
\par}
\label{f3c}
}
\hskip 0.4cm
\parbox{3in}{
\includegraphics[
height=2.0in,
width=2.8in
]%
{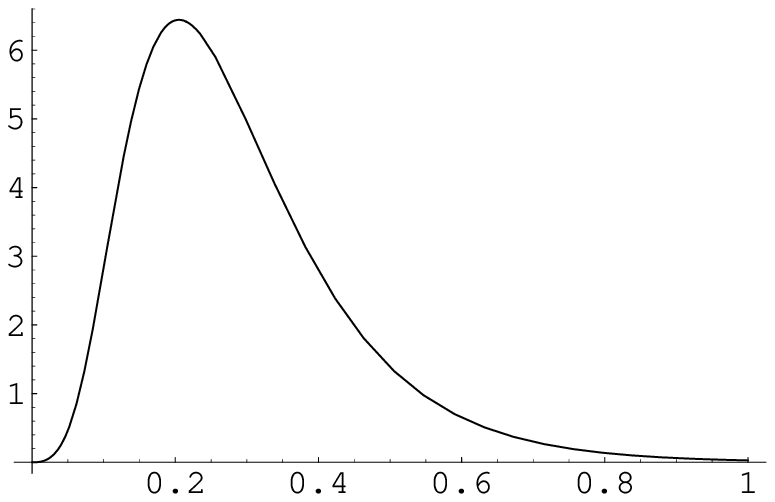}%
\par
{\footnotesize{
{\textbf{Figure 3f.~} $\nu=8$. Graph
of $\Rr_3(t)$, for $t \in [0,1)$.
One has $\Rr_3(0) = 0$, $\Rr_3(0.2) = 6.435...$,
$\Rr_3(0.4) = 2.787...$
$\Rr_3(0.6) = 0.6503...$,
$\Rr_3(0.9) = 0.06100...$.
The condition \rref{iprn3} for global existence of
the NS Cauchy problem,
i.e., $(\Dd_3 + \Rr_3)(t_1) \leqs \nu/G_3$,
is satisfied for any $t_1 \in [0.1567..., 1)$~.
}}
\par}
\label{f3f}
}
\end{figure}
\vfill\eject\noindent
\appendix
\section{Some comparison lemmas of the \v{C}aplygin type}
\label{appekapl}
Let $\Tl \in (0,+\infty]$; suppose we are given
\beq \ef \in C(\reali \times [0,\Tl), \reali),~(s, t) \mapsto \ef(s,t)~~
\mbox{such that}~~
{\partial \ef \over \partial s} \in C(\reali \times [0,\Tl), \reali)~; \feq
\beq s_0 \in \reali~. \feq
Under the above assumptions,
the Cauchy problem
\beq {d \Ss \over d t}(t) = \ef(\Ss(t),t)~, \qquad \Ss(0) = s_0~. \label{caus} \feq
has a unique maximal (i.e., nonextendable) solution $\Ss \in C^1([0,\Ts), \reali)$.
(Of course $\Ts \in (0,\Tl]$; later on we give conditions under which $\Ts=\Tl$.)
\par
The following is a known result, of the \v{C}aplygin type:
see \cite{Las} or \cite{Mitr}.
\begin{prop}
\label{lem1}
\textbf{Lemma.} Suppose there is a function $\Ww \in C([0,\Tl), \reali)$ such that
\beq {d_{+} \Ww \over d t}(t) \leqs \ef(\Ww(t),t)\quad \mbox{for $t \in [0,\Tl)$}, \qquad \Ww(0) \leqs s_0
\label{ipotesiw} \feq
(with  $d_{+}/d t$ the right lower Dini derivative, see Eq. \rref{dlower}). Then
\beq \Ww(t) \leqs \Ss(t) \quad \mbox{for $t \in [0,\Ts)$}~. \feq
\end{prop}
A straightforward variant of the previous result is the following.
\begin{prop}
\label{lem2}
\textbf{Lemma.} Suppose there is a function
$\Rr \in C([0,\Tl), \reali)$ such that
\beq {d^{+} \Rr \over d t}(t) \geqs \ef(\Rr(t),t)\quad \mbox{for $t \in [0,\Tl)$}, \qquad \Rr(0) \geqs s_0~,
\label{ipotesir} \feq
(with  $d^{+}/d t$ the right upper Dini derivative, see Eq. \rref{dupper}). Then
\beq \Rr(t) \geqs \Ss(t)~\mbox{for $t \in [0,\Ts)$}~. \label{so0} \feq
\end{prop}
\textbf{Proof.} We consider the function $-\Rr$, recalling that
${d_{+} (-\Rr)/d t} = - {d^{+} \Rr/d t}$; from \rref{ipotesir} we infer
\beq {d_{+} (-\Rr) \over d t}(t) \leqs -\ef(\Rr(t),t)~~\mbox{for $t \in [0,\Tl)$},
\qquad -\Rr(0) \leqs -s_0~. \feq
On the other hand, the function $-\Ss \in C^1([0,\Ts), \reali)$ is the maximal solution
of the Cauchy problem
\beq {d (-\Ss) \over d t}(t) = -\ef(\Ss(t),t)~, \qquad -\Ss(0) = -s_0~. \feq
Therefore, Lemma \ref{lem1} with $\Ww$ replaced by $-\Rr$ (and other obvious substitutions) gives
$$ -\Rr(t) \leqs -\Ss(t) \quad \mbox{for $t \in [0,\Ts)$}~, $$
yielding the thesis \rref{so0}. \fine
\begin{prop}
\label{lem3}
\textbf{Lemma.} Suppose there are functions $\Ww, \Rr \in C([0,\Tl), \reali)$ fulfilling
Eqs. \rref{ipotesiw}, \rref{ipotesir} respectively.
Then
\beq \Ts = \Tl~, \qquad \Ww(t) \leqs \Ss(t) \leqs \Rr(t) \quad\mbox{for $t \in [0,\Tl)$}~. \label{so} \feq
\end{prop}
\textbf{Proof.} From Lemmas \ref{lem1} and \ref{lem2} we infer
\beq \Ww(t) \leqs \Ss(t) \leqs \Rr(t) ~\mbox{for $t \in [0,\Ts)$}~; \label{so1} \feq
we claim that this inequality implies
\beq \Ts = \Tl~. \label{so2} \feq
In fact, whenever $\Ts$ is finite, the nonextendability  of $\Ss$
beyond this time implies that $\Ss$ is unbounded in any left
neighborhood of $\Ts^{-}$;
on the other hand, if it were $\Ts < \Tl$,
Eq. \rref{so1} would ensure the boundedness
of $\Ss$ on $[0, \Ts)$. \par
Now, the combination of \rref{so1} and \rref{so2} gives the thesis \rref{so}. \fine
If we forget the solution $\Ss$ of the Cauchy problem \rref{caus}, the inequality
\rref{so} becomes $\Ww(t) \leqs \Rr(t)$ for $t \in [0,\Tl)$; this is just the
statement in Eq. \rref{eqtesto} of Lemma \ref{lemtesto}, which is thus recognized as
a weaker version of Lemma \ref{lem3}.
\section{Appendix. Proof of Lemma
\ref{rem64} on the Galerkin approximants}
\label{appegal}
Let $\nu \in [0,+\infty)$, $f \in C([0,+\infty), \Dsz)$ and $u_0 \in \Dsz$.
In the finite dimensional space $\HG$, we consider
the maximal solution $\ug$ of the Cauchy problem \rref{galer}, here reproduced for
the reader's convenience:
$$ \mbox{Find $\ug \in C^1([0,\Tg), \HG)$ such that} $$
$$ {d \ug \over dt} = \nu \Delta \ug + \EG \PPP(\ug, \ug) +
\EG f~, \qquad \ug(0) = \EG u_0~. $$
For the moment, the maximal time of existence $\Tg$ is not known;
one of our aims is to show that $\Tg = +\infty$.
Let us proceed to prove this claim and all the other statements of
Lemma \ref{rem64}; our arguments are divided in several steps. \par \noindent
\textsl{Step 1 (``energy balance equation''). Everywhere on $[0, \Tg)$, one has}
\beq {d \over d t} \, \| \ug \|^2_{L^2} =
2 \nu \, \la \Delta \ug | \ug \ra_{L^2} +
2 \la \EG f | \ug \ra_{L^2}~. \label{ebal} \feq
In fact,
\beq {d \over d t}  \| \ug \|^2_{L^2} =
2 \la {d \ug \over d t} | \ug \ra_{L^2} \label{eb1} \feq
$$ = 2 \nu \la \Delta \ug | \ug \ra_{L^2} +
2 \la \EG \PPP(\ug, \ug) | \ug \ra_{L^2} + 2 \la \EG f | \ug \ra_{L^2}~. $$
On the other hand,
\beq \la \EG \PPP(\ug, \ug) | \ug \ra_{L^2} = \la \PPP(\ug, \ug)  | \ug \ra_{L^2} = 0~; \label{eb2} \feq
in the first passage above, the projection $\EG$ on $\HG$ is omitted
since we are taking the inner product with $\ug = \ug(t) \in \HG$;
the second passage reflects the well known identity $\la \PPP(v,w) | w \ra_{L^2} = 0$,
holding for all zero mean, divergence free vector fields $v, w$ on $\Td$
sufficiently regular to give meaning to the indicated
inner product (see, e.g., \cite{cog}, Lemma 2.3). Eqs. \rref{eb1} \rref{eb2}
yield the thesis \rref{ebal}.
\parn
\textsl{Step 2. Everywhere on $[0,\Tg)$, one has}
\beq {d \over d t} \, \| \ug \|^2_{L^2} \leqs
- 2 \nu \| \ug \|^2_{L^2} + 2   \| \EG f \|_{L^2} \| \ug \|_{L^2}~. \label{tps} \feq
This follows from the energy equation \rref{ebal},
combined with the inequalities
\beq \la \Delta \ug | \ug \ra_{L^2} \leqs - \| \ug \|^ 2_{L^2}~,
\qquad \la \EG f | \ug \ra_{L^2} \leqs  \| \EG f \|_{L^2} \| \ug \|_{L^2}~. \feq
The first of these relations follows using
\rref{inelap}, with $n=0$; the second one
is an application of the Cauchy-Schwarz inequality for $\la~|~\ra_{L^2}$.
\parn
\textsl{Step 3. Everywhere on $[0,\Tg)$, one has}
\beq {d^{+} \over d t} \, \| \ug \|_{L^2} \leqs
- \nu \| \ug \|_{L^2} + \| \EG f \|_{L^2}~, \label{tpr} \feq
\textsl{with $d^{+}/ d t$ the upper Dini derivative.}
Let us first prove Eq. \rref{tpr} at a time $t_0 \in [0, \Tg)$
such that $\ug(t_0) \neq 0$. In this case, $\| \ug(t) \|_{L^2}
\neq 0$ for all $t$ in an interval $I$ containing $t_0$;
in this interval the function $t \mapsto \| \ug \|_{L^2}$
is $C^1$, and we have
$$ {d^{+} \over d t} \, \| \ug \|_{L^2} =
{d \over d t} \, \| \ug \|_{L^2} = {1 \over 2 \| \ug \|_{L^2}}
{d \over d t} \, \| \ug \|^2_{L^2} \leqs
- \nu \| \ug \|_{L^2} + \| \EG f \|_{L^2}~, $$
the last passage following from \rref{tps}. Now,
consider an instant $t_0$ such that $\ug(t_0) =0$.
Due to a general result already mentioned
(see the comments before Eq. \rref{now1}), we can write
$$ \left. {d^{+} \over d t} \right|_{t_0} \, \| \ug \|_{L^2}
\leqs \| {d \ug \over d t}(t_0) \|_{L^2}~; $$
on the other hand, Eq. \rref{galer} and the assumption $\ug(t_0) = 0$
give  $(d \ug/d t)(t_0) = \EG f(t_0)$, so
$$ \left. {d^{+} \over d t} \right|_{t_0} \| \ug \|_{L^2}
\leqs \| \EG f(t_0) \|_{L^2}~; $$
this is just the thesis \rref{tpr}, at the instant under
consideration.
\par \noindent
\textsl{Step 4. One has}
\beq \| \ug(t) \|_{L^2} \leqs \Big(\| \EG u_0 \|_{L^2} + \int_{0}^t d s ~e^{\nu s}
\| \EG f(s) \|_{L^2} \Big) \, e^{-\nu t} \quad \mbox{for $t \in [0,\Tg)$}~. \label{jute} \feq
Let us consider the continuous function $t \in [0,\Tg) \mapsto \| \ug(t) \|_{L^2}$;
due to Step 3 and to $\ug(0) = \EG u_0$, we have
\beq {d^{+} \over d t} \| \ug(t) \|_{L^2} \leqs \ef(\|\ug(t) \|_{L^2}, t)~,
\qquad \| \ug(t) \|_{L^2} = \| \EG u_0 \|_{L^2} \feq
where
\beq \ef : \reali \times [0,+\infty) \vain \reali~, \qquad \ef(s,t) :=
- \nu s + \| \EG f(t) \|_{L^2}~. \feq
On the other hand, the Cauchy problem
\beq {d \Ss \over d t}(t) = \ef(\Ss(t), t)~, \qquad \Ss(0) = \| \EG u_0 \|_{L^2} \feq
has the global solution
\beq \Ss(t) := \Big(\| \EG u_0 \|_{L^2} + \int_{0}^t d s ~ e^{\nu s}
\| \EG f(s) \|_{L^2} \Big) \, e^{-\nu t} \quad \mbox{for $t \in [0,+\infty)$}~. \feq
So, by the comparison results reviewed
in the previous Appendix (see, in particular, Lemma \ref{lem2}), we have
$$ \| \ug(t) \|_{L^2} \leqs \Ss(t) \qquad \mbox{for $t \in [0, \Tg)$}~; $$
this is just the thesis \rref{jute}.
\par \noindent
\textsl{Step 5. One has $\Tg = +\infty$ (which justifies
statement \rref{tglob} in the lemma under proof).}
If it were $\Tg < + \infty$, we would have
$\limsup_{t \vain \Tg^{-}} \| \ug(t) \|_{L^2} = + \infty$;
this would be contradicted by Eq. \rref{jute}, implying
$\sup_{t \in [0, \Tg)} \| \ug(t) \|_{L^2} \leqs$
$\| \EG u_0 \|_{L^2} + \int_{0}^{\Tg} d s \, e^{\nu s}
\| \EG f(s) \|_{L^2}$ $< + \infty$.
\par \noindent
\textsl{Step 6. Proof of statement \rref{stacit} in the lemma,
under the assumption of zero forcing.} We use the results
of Steps 1 and 4, putting therein $\Tg = +\infty$.
For $f=0$ and $\nu=0$, the energy balance equation \rref{ebal} takes the form
$(d/d t) \, \| \ug \|^2_{L^2} = 0$; thus,
$$ \| \ug(t) \|_{L^2} = \| \ug(0) \|_{L^2} = \| \EG u_0 \|_{L^2} $$
for all $t \in [0,+\infty)$. For $f=0$ and $\nu > 0$, Eq.\rref{jute} gives
$$ \| \ug(t) \|_{L^2} \leqs \| \EG u_0 \|_{L^2} \, e^{-\nu t}~, $$
again for $t \in [0,+\infty)$. The above two equations correspond
to the content of \rref{stacit}. \parn
\textsl{Step 7. Proof of statements \rref{comb1}-\rref{comb3} in the lemma.}
Statement \rref{comb1} is just Eq. \rref{jute}, with $\Tg = + \infty$. Having proved \rref{comb1},
statements \rref{comb2} \rref{comb3} are obvious.
\vskip 0.7cm \noindent
\textbf{Acknowledgments.} This work was partly supported by INdAM, INFN and by MIUR, PRIN 2008
Research Project "Geometrical methods in the theory of nonlinear waves and applications".


\vskip 1cm \noindent
\begin{thebibliography}{99}
\bibitem{BKM} J. T.Beale, T. Kato, A. Majda, \textsl{Remarks on the breakdown of smooth solutions
for the 3D Euler equations}, Commun. Math. Phys. \textbf{94}, 61-66 (1984).
\vsm
\bibitem{Nec}
E. Behr, J. Ne$\check{\mbox{c}}$as, H. Wu,
\textsl{On blow-up of solution for Euler equations},
M2AN: Math. Model. Numer. Anal. \textbf{35} 229-238 (2001).
\vsm
\bibitem{Che} S.I. Chernyshenko, P. Constantin, J.C. Robinson, E.S. Titi,
\textsl{A posteriori regularity of the three-dimensional Navier-Stokes
equations from numerical computations}, J. Math. Phys. \textbf{48}, 065204/10 (2007).
\vsm
\bibitem{DR} M. Dashti, J.C. Robinson, \textsl{An a posteriori condition on the numerical
approximation of the Navier-Stokes equations for the existence of a strong solution},
SIAM J. Numer. Anal. \textbf{46}, 3136-3150 (2008).
\vsm
\bibitem{Kat3} H. Fujita, T. Kato, \textsl{On the Navier-Stokes initial value problem. I},
Arch. Rational Mech. Anal. \textbf{16}, 269--315 (1964).
\vsm
\bibitem{Kato} T.Kato, \textsl{Nonstationary flows of viscous and ideal fluids in $\reali^3$},
J.Funct.Anal. \textbf{9}, 296-305 (1972).
\vsm
\bibitem{Kat2} T. Kato, \textsl{Quasi-linear equations of evolution, with applications to
partial differential equations}, in ``Spectral theory and differential equations'',
Proceedings of the Dundee Symposium, Lecture Notes in Mathematics \textbf{448}, 23-70 (1975).
\vsm
\bibitem{Koz} H. Kozono, Y. Taniuchi, \textsl{Limiting case of the Sobolev inequality in BMO,
with application to the Euler equations}, Commun. Math. Phys. \textbf{214}, 191--200 (2000).
\bibitem{Las} V. Laksmikantham, S. Leela, ``Differential and
integral inequalities'', Volume I, Academic Press, New York (1969).
\vsm
\bibitem{Lem} P.G. Lemari\'e-Rieusset, ``Recent developments in the Navier-Stokes problem'', Chapman
\& Hall, Boca Raton (2002).
\vsm
\bibitem{Mitr} D.S. Mitrinovic, J.E. Pecaric, A.M. Fink, ``Inequalities
involving functions and their integrals and derivatives'', Kluwer, Dordrecht (1991).
\vsm
\bibitem{uno} C. Morosi, L. Pizzocchero, \textsl{On approximate solutions of semilinear
evolution equations}, Rev. Math. Phys. \textbf{16}, 383-420 (2004).
\vsm
\bibitem{due} C. Morosi, L. Pizzocchero, \textsl{On approximate solutions of
semilinear evolution equations II. Generalizations, and applications to Navier-Stokes
equations}, Rev. Math. Phys. \textbf{20}, 625-706 (2008).
\vsm
\bibitem{accau} C. Morosi, L. Pizzocchero,
\textsl{An $H^1$ setting for the Navier-Stokes equations: Quantitative estimates},
Nonlinear Anal. \textbf{74}, 2398-2414 (2011).
\vsm
\bibitem{cok} C. Morosi, L. Pizzocchero, \textsl{On the constants
in a basic inequality for the Euler and NS equations}, arXiv: 1007.4412v2 [mathAP] (2010).
\vsm
\bibitem{cog} C. Morosi, L. Pizzocchero, \textsl{On the constants
in a Kato inequality for the Euler and NS equations}, arXiv: 1009.2051v1 [mathAP] (2010).
To appear in Commun. Pure Appl. Anal.~.
\vsm
\bibitem{Pet} T. Petry, \textsl{On the stability of the Abramov transfer for differential-algebraic
equations of index 1}, SIAM J. Numer. Anal. \textbf{35}, 201--216 (1998).
\vsm
\bibitem{Rob} J.C. Robinson, W. Sadowski, \textsl{Numerical verification of
regularity in the three-dimensional Navier-Stokes equations for bounded sets of
initial data}, Asymptot. Anal. \textbf{59}, 39-50 (2008).
\vsm
\bibitem{Sin} Y. Sinai,
\textsl{Power series for solutions of the 3D Navier-Stokes
system on $\reali^3$}, J. Stat. Phys. \textbf{121}, 779--803 (2005).
\vsm
\bibitem{Tem} R. Temam, ``Navier-Stokes equations and nonlinear functional analysis'', SIAM,
Philadelphia (1983).
\end{thebibliography}
\end{document}